\documentclass[12pt]{amsart}
\usepackage{amssymb}

\ifx\mathrm\undefined\let\mathrm\rm\fi
\ifx\mathbf\undefined\let\mathbf\bf\fi
\ifx\mathfrak\undefined\let\mathfrak\frak\fi
\ifx\mathcal\undefined\let\mathcal\cal\fi
\ifx\mathbb\undefined\let\mathbb\Bbb\fi
\ifx\emph\undefined\let\emph\it\fi

\def\chgltxt#1{\marginpar{\kern-.5em$\blacklozenge$\kern.5em
\raggedright\footnotesize\bf#1}}

\newcommand{\g}{{{\mathfrak g}\,}}

\newcommand{\n}{{{\mathfrak n}}}
\newcommand{\h}{{{\mathfrak h\,}}}

\newcommand{\Z}{{\mathbb Z}}

\newcommand{\C}{{\mathbb C}}

\newcommand{\Ref}[1]{{(\ref{#1})}}

\newcommand{\la}{\lambda}

\newcommand{\dontprint}[1]

\newcommand{\nc}{\newcommand}

\newcommand{\bs}{\boldsymbol}
\nc{\Wr}{{ {\rm Wr}}}
\newcommand{\beq}{\begin{equation}}
\newcommand{\eeq}{\end{equation}}

\newcommand{\bean}{\begin{eqnarray}}
\newcommand{\eean}{\end{eqnarray}}
\newcommand{\be}{\begin{displaymath}}
\newcommand{\ee}{\end{displaymath}}
\newcommand{\bea}{\begin{eqnarray*}}
\newcommand{\eea}{\end{eqnarray*}}

{\relax}

\newtheorem%
{thm}{Theorem}[section]
\newtheorem%
{proposition}[thm]{Proposition}
\newtheorem%
{lemma}[thm]{Lemma}
\newtheorem%
{lemmadef}[thm]{Lemma-Definition}
\newtheorem%
{corollary}[thm]{Corollary}
\newtheorem%
{conjecture}[thm]{Conjecture}

\nc{\al}{\alpha}
\nc{\om}{\omega}
\nc{\La}{\Lambda}
\nc{\un}{U(\n_-)}
\nc{\px}{\frac d {d x}}

\newtheorem{theorem}{Theorem}[section]
\newtheorem{conj}[theorem]{Conjecture}
\newtheorem{cor}[theorem]{Corollary}
\newtheorem{lem}[theorem]{Lemma}
\newtheorem{prop}[theorem]{Proposition}

\theoremstyle{definition}

\theoremstyle{remark}

\theoremstyle{remark}

\numberwithin{equation}{section}

\nc{\on}{\operatorname}
\nc{\ch}{\mbox{ch}}
\nc{\pone}{{\mathbb C}{\mathbb P}^1}
\nc{\pa}{\partial}
\nc{\arr}{\rightarrow}
\nc{\larr}{\longrightarrow}
\nc{\ri}{\rangle}
\nc{\lef}{\langle}
\nc{\W}{{\mathcal W}}
\nc{\ep}{\epsilon}

\nc{\su}{\widehat{{\mathfrak sl}}_2}
\nc{\sw}{{\mathfrak s}{\mathfrak l}}

\nc{\G}{\widehat{\g}}
\nc{\De}{\Delta_+}
\nc{\gt}{\widetilde{\g}}
\nc{\Ga}{\Gamma}
\nc{\one}{{\mathbf 1}}
\nc{\z}{{\mathfrak Z}}
\nc{\Hh}{{\mathcal H}_\beta}
\nc{\qp}{q^{\frac{k}{2}}}
\nc{\qm}{q^{-\frac{k}{2}}}
\nc{\wt}{\widetilde}
\nc{\qn}{\frac{[m]_q^2}{[2m]_q}}
\nc{\cri}{_{\on{cr}}}
\nc{\kk}{h^\vee}
\nc{\sun}{\widehat{\sw}_N}
\nc{\hh}{\widehat{\mathfrak h}}
\nc{\HH}{{\mathcal H}_{q,t}}
\nc{\ca}{\wt{{\mathcal A}}_{h,k}(\sw_2)}
\nc{\gl}{\widehat{{\mathfrak g}{\mathfrak l}}_2}
\nc{\el}{\ell}
\nc{\s}{{\mathbf s}}
\nc{\bi}{\bibitem}
\nc{\WW}{\W_\beta}
\nc{\scr}{{\mathbf S}}
\nc{\ab}{{\mathbf a}}
\nc{\rr}{r}
\nc{\ol}{\overline}
\nc{\con}{qt^{-1} + q^{-1}t}
\nc{\den}{q^{\el-1} t^{-\el+1}+ q^{-\el+1} t^{\el-1}}
\nc{\ds}{\displaystyle}
\nc{\B}{B}
\nc{\A}{{\mathbb A}}
\nc{\GG}{{\mathcal G}}
\nc{\UU}{{\mathcal U}}
\nc{\MM}{{\mathcal M}}
\nc{\CC}{{\mathcal C}}
\nc{\dzz}{\frac{dz}{z}}
\nc{\Res}{\on{Res}}
\nc{\rep}{{\mathcal R}ep \;}
\nc{\uqg}{U_q \G}
\nc{\uqgg}{U_q \g}
\nc{\Fq}{{\mathbb F}_q}

\nc{\stimes}{\ltimes}
\nc{\K}{\hat{\mathcal K}}
\nc{\Ql}{\ol{\mathbb Q}_\ell}

\nc{\ga}{\gamma}
\nc{\PL}{{}^L P}
\nc{\E}{\mc E}
\nc{\mc}{\mathcal}
\nc{\mbf}{\mathbf}
\nc{\OO}{{\mc O}}
\nc{\Po}{{\mc P}}
\nc{\V}{{\mc V}}
\nc{\yy}{{\mc Y}}
\nc{\M}{\mathcal M}
\nc{\Coh}{{{\mathcal C}oh}}
\nc{\Cohn}{\Coh_n}
\nc{\f}{{\mathcal F}}
\nc{\si}{_E}
\nc{\Gaf}{{\mathbb G}_{a,\Fq}}
\nc{\KK}{{\mathfrak k}}

\nc{\PCr}{{ \bs P (\C[x])^r }}
\nc{\PCN}{{ \bs P (\C[x])^N }}

\nc{\sN}{sl_{2N+1}}
\nc{\Pzr}{{ \bs P(\C((x-z)))^r}}
\nc{\PzN}{{ \bs P(\C((x-z)))^N}}

\newcommand{\p}{\partial_x}
\newcommand{\Om}{\Omega}
\newcommand{\ox}{\otimes}

\def\Vb{L_\bullet}
\def\gl{\frak{gl}}
\def\glt{\gl_2}
\def\glN{\gl_N}
\def\glM{\gl_M}
\def\slN {\frak{sl}_N}
\def\slM {\frak{sl}_M}
\def\V-{\>\hbox{$\=V\}$-}}

\def\fratop{\genfrac{}{}{0pt}1}

\newcommand{\NNN}{^{\langle N\rangle}}
\newcommand{\MMM}{^{\langle M\rangle}}
\newcommand{\Zp}{\Z_{\geq 0}}

\newcommand{\GR}{{\rm Gr}(X,N)}
\newcommand{\GRM}{{\rm Gr}(Y,M)}

\newcommand{\Wrd}{{\rm Wr}^{\rm(d)}}

\newcommand{\LMN}{\bs L_{\bs m}[\bs n]}
\newcommand{\LNM}{\bs L_{\bs n}[\bs m]}

\newcommand{\Sym}{{\rm Sym}}

\newcommand{\Gr}{{\rm Gr}}

\newcommand{\tn}{{\bs t^{\langle \bs n\rangle}}}

\begin{document}

\title[Bispectral and $(\glN,\glM)$ Dualities]
{Bispectral and $(\glN,\glM)$ Dualities,\\ Discrete versus Differential}

\author[E. Mukhin, V. Tarasov, and A. Varchenko]
{E. Mukhin ${}^{*}$, V. Tarasov ${}^{*,\star,1}$,
\and A. Varchenko {${}^{**,2}$} }
\thanks{${}^1$\ Supported in part by RFFI grant 05-01-00922}
\thanks{${}^2$\ Supported in part by NSF grant DMS-0244579}

\maketitle

\centerline{\it ${}^*$Department of Mathematical Sciences,
Indiana University -- Purdue University,}
\centerline{\it Indianapolis, 402 North Blackford St, Indianapolis,
IN 46202-3216, USA}
\smallskip
\centerline{\it $^\star$St.\,Petersburg Branch of Steklov Mathematical
Institute}
\centerline{\it Fontanka 27, St.\,Petersburg, 191023, Russia}
\smallskip
\centerline{\it ${}^{**}$Department of Mathematics, University of
North Carolina at Chapel Hill,} \centerline{\it Chapel Hill, NC
27599-3250, USA} \medskip

\medskip
\begin{center}
May, 2006
\end{center}

\thispagestyle{empty}

\begin{abstract}
Let $V = \langle\, x^{\la_i}p_{ij}(x),\ i=1,\dots,n, \ j=1, \dots ,
N_i\, \rangle$ be a space of quasi-polynomials in $x$ of dimension
$N=N_1+\dots+N_n$. The regularized fundamental differential operator of $V$
is the polynomial differential operator $\sum_{i=0}^N
A_{N-i}(x)(x\p)^i$ annihilating $V$ and such that its leading
coefficient $A_0$ is a monic polynomial of the minimal possible degree. Let
$U = \langle \,z_a^{u}\,q_{ab}(u),\, a=1,\dots,m,\ b=1,\dots ,M_a
\rangle$ be a space of quasi-exponentials in $u$ of dimension
$M=M_1+\dots+M_m$. The regularized fundamental difference operator of $U$
is the polynomial difference operator $\sum_{i=0}^M
B_{M-i}(u)(\tau_u)^i$ annihilating $U$ and such that its leading
coefficient $B_0$ is a monic polynomial of the minimal possible degree.
Here $(\tau_uf)(u)=f(u+1)$.

Having a space $V$ of quasi-polynomials with the regularized
fundamental differential operator $D$, we construct a space of
quasi-exponentials $U = \langle \,z_a^{u}\,q_{ab}(u)\, \rangle$ whose
regularized fundamental difference operator is the difference operator
$\sum_{i=0}^N u^i A_{N-i}(\tau_u)$. The
space $U$ is constructed from $V$ by a suitable integral transform.
Similarly, having $U$ we can recover $V$ by a suitable integral transform.
Our integral transforms are analogs of the bispectral involution on
the space of rational solutions to the KP hierarchy \cite{W}.

As a corollary of the properties of the integral transforms
we obtain a correspondence between solutions to the Bethe ansatz equations
of two $(\glN,\glM)$ dual quantum integrable models: one is the
special trigonometric Gaudin model and the other is the special XXX model.

\end{abstract}

\section{Introduction}
Let $V = \langle\, x^{\la_i}p_{ij}(x),\ i=1,\dots,n, \ j=1, \dots ,
N_i\, \rangle$ be a space of quasi-polynomials in $x$ of dimension
$N=N_1+\dots+N_n$. The regularized fundamental differential operator of $V$
is the polynomial differential operator
$\sum_{i=0}^N A_{N-i}(x)(x\p)^i$
annihilating $V$ and such that its leading
coefficient $A_0$ is a monic polynomial of the minimal possible degree.

Let
$U = \langle \,z_a^{u}\,q_{ab}(u),\, a=1,\dots,m,\ b=1,\dots ,M_a
\rangle$ be a space of quasi-exponentials in $u$ of dimension
$M=M_1+\dots+M_m$. The regularized fundamental difference operator of $U$
is the polynomial difference operator
$\sum_{i=0}^M B_{M-i}(u)(\tau_u)^i$
annihilating $U$ and such that its leading
coefficient $B_0$ is a monic polynomial of the minimal possible degree.
Here $(\tau_uf)(u)=f(u+1)$.

Having a space $V$ of quasi-polynomials with the regularized
fundamental differential operator
$\sum_{i=0}^N A_{N-i}(x)(x\p)^i$,
we construct a space of
quasi-exponentials $U = \langle \,z_a^{u}\,q_{ab}(u)\, \rangle$ whose
regularized fundamental difference operator is the difference operator
\linebreak
$\sum_{i=0}^N u^i A_{N-i}(\tau_u)$. The
space $U$ is constructed from $V$ by a suitable integral transform.

Having a space $U = \langle \,z_a^{u}\,q_{ab}(u)\, \rangle$
of quasi-exponentials with the regularized
fundamental difference operator
$\sum_{i=0}^M B_{M-i}(u)(\tau_u)^i$,
we construct a space of
quasi-polynomials $V = \langle\, x^{\la_i}p_{ij}(x) \rangle$
whose regularized fundamental differential operator is the differential operator
$\sum_{i=0}^M x^i B_{M-i}(x\p)$.
The space $V$ is constructed from $U$ by a suitable integral transform.

Our integral transforms are analogs of the bispectral involution on
the space of rational solutions to the KP
hierarchy \cite{W}.

As a corollary of the properties of our integral transforms
we obtain a correspondence between solutions to the Bethe ansatz equations
of two $(\glN,\glM)$ dual quantum integrable models: one is the
special trigonometric Gaudin model and the other is the special XXX model.

\bigskip

\noindent
{\bf Example.}
Let $\bs n=(n_1,n_2)$ and $\bs m=(m_1,m_2)$ be two vectors of nonnegative integers such that
$n_1+n_2=m_1+m_2$.
Let $d$ be the number of integers $i$ such that
$\max\,(0, n_2-m_1) \leq i\leq \min \,(m_2,n_2)$.
Let $\bs z = (z_1,z_2), \, \bs \la = (\la_1, \la_2) \,\in\,
\C^2$\ be two points with distinct coordinates.

Consider two systems of algebraic equations. The first system is the system of equations
\bean
\label{bea N_2}
&&
\frac{\la_1-\la_2-1}{t_i}\ +\
\sum_{a=1}^2
\frac{m_a}{t_i-z_a} \
-\ \sum_{j=1,\ j\neq i}^{n_2}
\frac{2}{t_i-t_{j}}\ =\ 0\ ,
\quad
i=1,\dots, n_2\ ,
\eean
with respect to the
unknown numbers
$t_1,\dots,t_{n_2}$. The system is symmetric with respect to the group
$\Sigma_{n_2}$ of permutations of $t_1,\dots,t_{n_2}$.
One can show that for generic $\bs z, \bs \la$
the number of $\Sigma_{n_2}$-orbits of solutions of system \Ref{bea N_2} is equal to $d$.
This system is called the system of the
Gaudin Bethe ansatz equations, see Section \ref{sec crit points}.

The second systems of equations is the system
\bean
\label{bea M_2}
&&
\prod_{i=1}^2 \ \frac {s_a -\la_i - 1}
{s_a-\la_i - 1 - n_i } \
\prod_{b, \ b \neq a}^{m_2}\
\frac {s_a - s_{b} - 1}
{s_b - s_{b} + 1} \
= \ \frac{z_2}{z_1}\ ,
\qquad
a=1,\dots,m_2\ ,
\eean
with respect to the unknown numbers
$s_1,\dots,s_{m_2}$. The system is symmetric with respect to the group
$\Sigma_{m_2}$ of permutations of $s_1,\dots,s_{m_2}$.
One can show that for generic $\bs z, \bs \la$
the number of $\Sigma_{m_2}$-orbits of solutions of system \Ref{bea M_2} is equal to $d$.
This system is called the system of the
XXX Bethe ansatz equations, see Section \ref{sec crit points}.

\medskip
\noindent
{\bf Problem.} Establish a correspondence between orbits of solutions to systems
\Ref{bea N_2} and \Ref{bea M_2}.

We give two solutions to this problem.

\medskip
\noindent
{\it The first solution.} To system \Ref{bea N_2}, we assign the vector space
$(L_{m_1}\ox L_{m_2})[n_1,n_2]$ and four commuting linear operators
acting on this space.
Here $(L_{m_1}\ox L_{m_2})[n_1,n_2]$ denotes the weight subspace of weight
$[n_1,n_2]$ of the tensor product of $\glt$-modules with highest weights
$(m_1,0)$ and $(m_2,0)$, respectively. The space $(L_{m_1}\ox L_{m_2})[n_1,n_2]$
is of dimension $d$.
The linear operators are denoted by
\bea
H_1^{\frak G}(\la_1+n_1, \la_2+n_2, \bs z)\ ,
&\qquad &
H_2^{\frak G}(\la_1+n_1, \la_2+n_2, \bs z)\ ,
\\
G_1^{\frak G}(\la_1+n_1, \la_2+n_2, \bs z)\ ,
&\qquad &
G_2^{\frak G}(\la_1+n_1, \la_2+n_2, \bs z)
\eea
and called the Gaudin
KZ and dynamical Hamiltonians. To each orbit of solutions of system \Ref{bea N_2}, the Bethe
ansatz method assigns an eigenvector of the commuting Hamiltonians. The constructed
Bethe vectors form a basis of this vector space.

To system \Ref{bea M_2}, we assign the vector space
$(L_{n_1}\ox L_{n_2})[m_1,m_2]$ of the same dimension $d$
and four commuting linear operators
acting on this space.
The linear operators are
denoted by
\bea
H_1^{\frak X}(\bs z, \la_1+n_1, \la_2+n_2)\ ,
&\qquad &
H_2^{\frak X}(\bs z, \la_1+n_1, \la_2+n_2)\ ,
\\
G_1^{\frak X}(\bs z, \la_1+n_1, \la_2+n_2)\ ,
&\qquad &
G_2^{\frak X}(\bs z, \la_1+n_1, \la_2+n_2)
\eea
and
called the XXX KZ and dynamical Hamiltonians.
To each orbit of solutions of system \Ref{bea M_2}, the Bethe
ansatz method assigns an eigenvector of the commuting Hamiltonians. The constructed
Bethe vectors form a basis of this vector space.

There is a natural isomorphism of the vector spaces $(L_{m_1}\ox L_{m_2})[n_1,n_2]$
and $(L_{n_1}\ox L_{n_2})[m_1,m_2]$, which identifies the Hamiltonians:
\bea
H_a^{\frak G}(\la_1+n_1, \la_2+n_2, \bs z)\ &=&\
G_a^{\frak X}(\bs z, \la_1+n_1, \la_2+n_2)\ ,
\\
G_i^{\frak G}(\la_1+n_1, \la_2+n_2, \bs z)\ &=&\
H_i^{\frak X}(\bs z, \la_1+n_1, \la_2+n_2)\ ,
\eea
for $i, a = 1, 2$. This isomorphism is called the $(\glt,\glt)$-duality, see \cite{TV5}.
Under the duality isomorphism the Bethe vectors are identified and this
identification gives a correspondence
between the orbits of solutions of systems \Ref{bea N_2} and \Ref{bea M_2}.

\medskip
\noindent
{\it The second solution.} To the orbit of a solution $(t_1,\dots,t_{n_2})$ of
system \Ref{bea N_2}, we assign the polynomial $p_2(x) = \prod_{i=1}^{n_2} (x-t_i)$
and the differential operator
\bea
D = x^2 (x-z_1)(x-z_2)
( \partial_x -
\ln' ( \frac { x^{\la_1-1}\prod_{a=1}^2 (x-z_a)^{m_a} } { p_2 } ) )
\left( \partial_x - \ln' ( x^{\la_2} p_2 ) \right) \ ,
\eea
where $\p = d/dx$ and
for any function $f$,\ $\ln' f$ denotes the logarithmic derivative $f'/f$.
Clearly the differential equation
$ D f(x) = 0$ has a solution $x^{\la_2}\, p_2(x)$.

We show that $D$ can be written in the form
\bea
D\ =\
A_0(x) (x\p)^2 + A_1(x) x\p + A_2(x)
\eea
where $A_0, A_1, A_2$ are
polynomials in $x$ of degree not greater than two. Then we consider the second order
difference equation
\bea
(\,u^2\, A_0(\tau_u) \,+ \, u\,A_1(\tau_u)\, + A_2(\tau_u)\,)\, g(u)\ =\ 0\
\eea
for the unknown function $g(u)$.
It turns out that this difference equation
has a solution of the form $z_2^u\,q_2(u)$, where
$q_2(u) = \prod_{a=1}^{m_2}\, (u-s_a)$, and the roots $s_1,\dots,s_{m_2}$
form a solution of system \Ref{bea M_2}. This construction gives a second correspondence
between the orbits of solutions of systems \Ref{bea N_2} and \Ref{bea M_2}.

\medskip

We show that the two described constructions (the first, based on the
$(\glt,\glt)$-duality, and the second, which uses the differential and
difference operators) give the same correspondence between solutions of systems
\Ref{bea N_2} and \Ref{bea M_2}.

\bigskip

This paper is a development of results of the paper \cite{MTV1}, in which
we presented an integral transform giving an involution on the space of
quasi-exponentials, the involution which corresponds to the bispectral
involution of G.\,Wilson in \cite {W}.

\bigskip

The paper has the following structure. In Section \ref{sec quasi-polynomials},
we discuss spaces of quasi-polynomials and their
fundamental differential operators. In Section \ref{sec
quasi-exponentials}, we discuss spaces of quasi-exponentials and their
fundamental difference operators. In Section \ref{Sec Int
transforms}, we define integral transforms establishing a duality
between spaces of quasi-polynomials and quasi-exponentials. In
Section \ref{special spaces}, we introduce special spaces of
quasi-polynomials and quasi-exponentials. We introduce special
integral transforms relating the special spaces. In Section
\ref{Special spaces and solutions of the Bethe ansatz equations}, we
discuss relations between special spaces and solutions of the Bethe
ansatz equations. The special spaces of quasi-polynomials correspond
to solutions of the Bethe ansatz equations in (trigonometric) Gaudin
models. The special spaces of quasi-exponentials correspond to
solutions of the Bethe ansatz equations in XXX models. The number of
solutions of the Bethe ansatz equations is discussed in Section
\ref{Sec FINIT}. We describe the Gaudin and XXX models in Section
\ref{section KZ hamiltonians}. The $(\glN, \glM)$ duality between the
(trigonometric) Gaudin and XXX models is defined in Section \ref{sec
duality}. In Section \ref{sec duality}, we formulate a conjecture
about the bispectral correspondence of Bethe vectors under the $(\glN,
\glM)$ duality. In Section \ref{sec m=n}, we prove the conjecture for
$N=M=2$. In Section \ref{sec baker}, we discuss the relation of our
integral transforms to the bispectral correspondence of suitable
Baker-Akhieser functions for differential and difference operators.

\medskip

Authors thank Y.\,Berest and O.\,Chalykh for discussions.

\section{Spaces of quasi-polynomials }
\label{sec quasi-polynomials}

\subsection{Definition}
\label{subs quasi-polynomials}
Let $p\in\C[x]$ be a polynomial, $\la$ a complex number.
The function $x^\la p$\ is called {\it a quasi-polynomial in $x$}. The
quasi-polynomial is a multi-valued function. Different local
uni-valued branches of the function differ by a nonzero constant
factor.

Let $N_1, \dots , N_n$ be natural numbers. Set $N = N_1 + \dots +
N_n$. For $i = 1,\dots, n$,\ let \ $0 < n_{i1} < \dots < n_{iN_i}$\
be a sequence of positive integers. For $i= 1,\dots , n,\ j = 1,
\dots , N_i$, \ let\ $p_{ij}\in \C[x]$\ be a polynomial of degree
$n_{ij}$.

Let $\la_1, \dots , \la_n\,\in\, \C$ be distinct
numbers such that $\la_i-\la_j \notin \Z$ for $i\neq j$.

Denote by $V$ the complex vector space spanned by functions \
$x^{\la_i}p_{ij}(x)$,\ $i=1,\dots , n$,\linebreak $j=1,\dots ,
N_i$. The dimension of $V$ is $N$.

The space $V$ is called {\it the space of quasi-polynomials.}

We say that the space $V$ is {\it non-degenerate} if
for any $i=1,\dots , n$ and any $j=1,\dots , N_i$,
\begin{enumerate}
\item[$\bullet$]
there exists
a linear combination of polynomials $p_{i1},p_{i2},\dots,p_{iN_i}$ which has a root
at $x=0$ of multiplicity $j-1$,
\item[$\bullet$]
the space $V$ does not contain the function $x^{\la_i +n_{ij}}$.
\end{enumerate}

\subsection{Exponents}
\label{subsec exponents}
Let $V$ be the space of quasi-polynomials.

For $z\in \C^*$, define {\it the sequence of exponents of $V$ at $z$}
as the unique sequence of integers, $\bs e = \{e_1 < \dots < e_N\}$,
with the property: for $i=1,\dots,N$, there exists $f\in V$ such that
$f$ has a root at $x=z$ of multiplicity $e_i$.

We say that $z\in \C^*$ is a singular point of $V$ if the set of
exponents of $V$ at $z$ differs from the set $\{0,\dots, N-1\}$. The
space $V$ has at most finitely many singular points.

Let $(z_1,\dots ,z_m)$ be the subset of $\C^*$ of all singular points of $V$.
For $a=1,\dots,m$, let
$$
\{\ 0 < \dots < N - M_a -1 < N - M_a + m_{a 1} < \dots < N - M_a + m_{a M_a}\ \}
$$
be the exponents of $V$ at $z_a$. Here $0<m_{a 1}<\dots < m_{a
M_a}$ and $M_a$ is an integer such that $1\leq M_a \leq N$. Set $M =
M_1 + \dots + M_m$.

We say that $V$ is {\it a space of quasi-polynomials with data}
\bea
\frak D_V\ =\ \{\,n, N_i,n_{ij},\la_i,\,m, M_a, m_{ab}, z_a\,\}
\eea
where $i=1,\dots, n,\, j=1,\dots,N_i,\, a=1,\dots,m,\, b=1,\dots,M_a$.

\subsection{Fundamental differential operator}
\label{subsection fund oper}
For functions $f_1,\dots, f_i$ of one variable, denote by
$\Wr\,(f_1,\dots,f_i)$ their Wronskian, that is, the determinant of the
$i\times i$-matrix whose $j$-th row is $f_j, f_j^{(1)},\dots,f_j^{(i-1)}$.

Define the Wronskian of $V$, denoted by $\Wr_V$, as the Wronskian of a
basis of $V$. The Wronskian of $V$ is determined up to multiplication
by a nonzero number.

\begin{lem}\label{Wr. identity}
Let $V$ be a non-degenerate space of quasi-polynomials, then
\bea
\sum_{a=1}^m\sum_{b=1}^{M_a}\ (\,m_{a b} +1-b\,)
\ = \ \sum_{i=1}^n\sum_{j=1}^{N_i} \ (\,n_{i j}+1-j\,) \ .
\eea
\end{lem}
The lemma is proved by analyzing the order of zeroes of the Wronskian of $V$
and its asymptotics at infinity.

\bigskip

{\it The monic fundamental differential operator of $V$} is the unique monic
linear differential operator of order $N$ annihilating $V$. It is denoted by
$\bar D_V$. We have
\bea
\bar D_V\ =\ \p^N + \bar
A_1\p^{N-1} + \dots + \bar A_N\ , \qquad \bar A_i\ =\ (-1)^i\, \frac
{\Wr_{V, i}}{\Wr_V}\ ,
\eea
where $\p = d/dx$,\ {} $\Wr_V$ is the
Wronskian of a basis $\{f_1,\dots ,f_N\}$ of $ V$, \ {} $\Wr_{V,i}$ is the
determinant of the $N\times N$-matrix whose
$j$-th row is $f_j,
f_j^{(1)},\dots,f_j^{(N-i-1)},f_j^{(N-i+1)},\dots , f_j^{(N)}$.

For any $j=1,\dots,N$, the order of the pole of $\bar A_j$ at $x=z_a$,
$a=1,\dots,m$, does not exceed $j$.

\begin{lem}\label{irreg of fund oper}
Let $V$ be a non-degenerate space of quasi-polynomials, then
for $a=1,\dots,m$, the order of the pole of
$\bar A_{ M_a}$ at $x=z_a$ is $ M_a$ and the order of the pole
of $\bar A_i$ is not greater than $ M_a$ for $i > \bar M_a$.
\end{lem}
The proof follows from counting orders of determinants $W_{V,i}$.

The polynomial differential operator
\bea
D_V\ = \ x^N \prod_{a=1}^m\ (x - z_a)^{M_a}\
\bar D_V
\eea
is called {\it the regularized fundamental differential operator of $V$}.

It is easy to see that the regularized fundamental differential operator $D_V$ of
$V$ has the form $\tilde A_0(x)x^N \p^N\!+\tilde
A_1(x)x^{N-1}\p^{N-1}\!+\dots +\tilde A_N(x)$ where $\tilde A_i(x)$ are
polynomials. Using the formula $x^i\p^i = x\p( x\p-1)\cdots (x\p-i+1)$
we can present the regularized fundamental differential operator in the form
\bea
D_V = A_0(x) (x\p)^N\!+A_1(x)(x\p)^{N-1}\!+ A_2(x)(x\p)^{N-2}\!+\dots +A_N(x)
\eea
with polynomial coefficients $A_i(x)$, $i = 0,\dots,N$.

\begin{lem}\label{deg coeff}
Let $V$ be a non-degenerate space of quasi-polynomials, then

\begin{enumerate}
\item[$\bullet$]
We have
$
A_0(x) \ =\ \prod_{a=1}^m\ (x-z_a)^{M_a}.
$
\item[$\bullet$] All coefficients $A_i(x)$ are polynomials in $x$
of degree not greater than $M$.

\item[$\bullet$]
Write $D_V = x^M B_0(x\p) + x^{M-1} B_1(x\p) + \dots + x B_{M-1}(x\p)
+ B_M(x\p)$ where $B_0(x\p),\dots, B_M(x\p)$ are polynomials in $x\p$ with constant
coefficients. Then
\bea
B_0(x\p) = \prod_{i=1}^n\prod_{j=1}^{N_i} (x\p - \la_i-n_{ij})\ ,
\quad
B_M(x\p) = \prod_{a=1}^m(-z_a)^{M_a}
\prod_{i=1}^n\prod_{j=1}^{N_i} (x\p - \la_i-j+1)\ .
\eea
\item[$\bullet$] The polynomials $B_0,\dots, B_M$ have no
common factor of positive degree.

\end{enumerate}
\hfill
$\square$

\end{lem}

\subsection{Conjugate Space}\label{subsection conj space}
Let $V$ be a space of quasi-polynomials as in Section \ref{subs quasi-polynomials}.

The complex vector space spanned by all functions of the form
$\Wr(f_1,\dots,f_{N-1})/\Wr_V$ with $f_i\in V$ has dimension $N$. It is denoted by
$V^\star$ and called {\it conjugate to} $V$.

The complex vector space spanned by all functions of the form
\linebreak
$f(x)\, x^{-N} \prod_{a=1}^m(x-z_a)^{-M_a}$\ with $f\in V^\star$
is denoted by
$V^\dagger$ and called {\it regularized conjugate to} $V$.

\begin{lem}\label{conj exponents}
Let $V$ be a non-degenerate space of quasi-polynomials, then
\begin{enumerate}
\item[(i)]
For $a=1, \dots , m$, if
$\bs e$ are exponents of $V$ at $z_a$, then
\bea
\bs e^\star = \{ - e_{N}-1+N < - e_{N-1}-1+N < \dots < - e_{1}-1+N\}
\eea
are exponents of $V^\star$ at $z_a$ and
\bea
\bs e^\dagger = \{ - e_{N}-1+N-M_a < - e_{N-1}-1+N-M_a < \dots < - e_{1}-1+N-M_a\}
\eea
are exponents of $V^\dagger$ at $z_a$.
\item[(ii)]
For any $i=1,\dots,n$, $j=1,\dots,N_i$, there exists $f\in V^\dagger$ such
that the function $x^{\la_i+j}f$ has a nonzero limit as $x\to 0$.
\item[(iii)]
For any $i=1,\dots,n$, $j=1,\dots,N_i$, there exists $f\in V^\dagger$ such
that the function $x^{\la_i+n_{ij}+1-N_i}f$ has a nonzero limit as $x\to \infty$.
\end{enumerate}

\hfill
$\square$
\end{lem}

Let $D = \sum_i A_i(x) \p^i$ be a differential operator with meromorphic coefficients.
The operator $D^* = \sum_i(-\p)^i A_i(x)$ is called {\it
formal conjugate to} $D$.

\begin{lem}\label{conj mon oper}
Let $V$ be a non-degenerate space of quasi-polynomials.
Let $\bar D_V$ and $D_V$ be the monic
and regularized fundamental differential operators of $V$, respectively. Then
$(\bar D_V)^*$ annihilates $V^\star$ and $(D_V)^*$ annihilates $V^\dagger$.
\hfill
$\square$
\end{lem}

\section{Spaces of quasi-exponentials }
\label{sec quasi-exponentials}

\subsection{Definition}
\label{subsec quasi-exponentials}
Define the operator $\tau_u$ acting on functions of $u$ as
$(\tau_uf)(u)=f(u+1)$.

Let $z$ be a nonzero complex number with fixed argument.
Set $z^u = e^{u\ln\,z}$. We have $\tau_u z^u = z^uz$.

Let $q\in\C[u]$ be a polynomial.
The function $z^u q $ is called {\it a quasi-exponential in $u$}.

\bigskip

Let $M_1, \dots , M_m$ be natural numbers. Set $M = M_1 + \dots +
M_m$. For $a = 1,\dots, m$,\ let \ $0 < m_{a1} < \dots < m_{aM_a}$\
be a sequence of positive integers. For $a= 1,\dots , m,\ b = 1,
\dots , M_a$, \ let\ $q_{ab}\in \C[u]$\ be a polynomial of degree
$m_{ab}$.

Let $z_1, \dots , z_m$ be distinct nonzero complex numbers with fixed
arguments.

Denote by $U$ the complex vector space spanned by functions \ $z_a^{u}q_{ab}(u) $,\
$a=1,\dots , m$,\linebreak $b=1,\dots , M_a$. The dimension of $U$ is $M$.

The space $U$ is called {\it the space of quasi-exponentials.}

\bigskip

For functions $f_1,\dots, f_i$ of $u$, denote by $\Wrd(f_1,\dots,f_i)$
their discrete Wronskian which is the determinant of the $i\times i$-matrix
whose $j$-th row is $f_j(u), f_j(u+1),\dots,f_j(u+i-1)$.

Define {\it the discrete Wronskian of $U$}, denoted by $\Wrd_U$, as the
discrete Wronskian of a basis of $U$. The discrete Wronskian of $U$
is determined up to multiplication by a nonzero constant.

\begin{lem}
\label{lem on wronskian}
We have
$$
\Wrd_U(u)\  =\ S(u)\, \prod_{a=1}^m \,z_a^{M_au}\ ,
$$
where $S(u)$ is a polynomial of degree
\ {} $\sum_{i=1}^n\sum_{j=1}^{N_i}(n_{ij}+1-j)$.
\hfill $\square$
\end{lem}

\subsection{The frame of a space of quasi-exponentials}
Let $U$ be a space of quasi-expo-
nentials like in Section \ref{subsec quasi-exponentials}.
Let $v_1,\dots,v_M$ be the quasi-exponentials
generating $U$.

For $i=1,\dots,M$, let $S_i\in \C[u]$ be the monic polynomial of the
greatest degree such that the function $\Wrd(v_{j_1},v_{j_2},\dots,v_{j_i})/S_i$ is
regular for any $j_1,j_2,\dots,j_i\in\{1,\dots,M\}$.

In particular, for the discrete Wronskian of $U$ we have
\bea
\Wrd_U(u)\ =\ {\rm const}\ S_M(u)\ \prod_{a=1}^m\, z_a^{M_au}
\eea
with a nonzero constant.

\begin{lem}
\label{U=T}
There exists a unique sequence of monic polynomials
$P_1(u),\dots,P_M(u)$ such that
\be
S_i(u)\ =\ \prod_{k=1}^i \prod _{j=1}^{i-k+1}P_k(u+j-1)
\ee
for $i=1,\dots,M$.
\end{lem}

This lemma is an analog of Lemma 4.12 in \cite{MV2}.

The monic polynomials $P_1(u), \dots , P_M(u)$ are called {\it the
frame} of $U$.

\begin{proof}
We construct $P_i$ by induction on $i$. For $i=1$, we set $P_1=S_1$.
Suppose the lemma is proved for all $i=1,\dots, i_0-1$. Then we set
\be
R(u)\ =\ \prod_{i=1}^{i_0-1}\prod _{j=1}^{i_0-i+1}P_i(u+j-1)\ ,
\qquad
P_{i_0}(u)\ =\ S_{i_0}(u)/R(u)\ .
\ee
We have to show that
$P_{i_0}$ is a polynomial. In other words, we have to show that the
discrete Wronskian of any $i_0$-dimensional subspace $\langle
v_{j_1},\dots,v_{j_{i_0}}\rangle$ is divisible by $R(u)$.

Consider the Grassmannian $Gr(i_0-2,U)$ of $(i_0-2)$-dimensional spaces
in $U$. For any $z\in \C$ the set of points in $Gr(i_0-2,U)$ such that
the corresponding discrete Wronskian divided by $S_{i_0-2}$ does not vanish at
$z$, is an open set. Therefore, we have an open
set of points in $Gr(i_0-2,U)$ such that the corresponding discrete
Wronskian divided by $S_{i_0-2}$ does not vanish at roots of $R(u-1)$. We call
such subspaces acceptable.

Therefore, we have an open set of points in $Gr(i_0,U)$ such
that the corresponding $i_0$-dimensional space contains an acceptable
$i_0-2$ dimensional subspace. Let $w_1,\dots,w_{i_0}\in U$ be such
that $w_1,\dots,w_{i_0-2}$ span an acceptable space.
It is enough to show that $\Wrd(w_1,\dots,w_{i_0})$ is divisible by
$R(u)$.

Using discrete Wronskian identities of \cite{MV2},
we have for suitable holomorphic
functions $f_1,f_2,g$:
\bea
\Wrd(w_1,\dots,w_{i_0})=\frac{\Wrd(\Wrd(w_1,\dots,w_{i_0-1}),
\Wrd(w_1,\dots,w_{i_0-2},w_{i_0}))}{\Wrd(w_1,\dots,w_{i_0-2})(u+1)} =
\\
\frac{\Wrd(S_{i_0-1}f_1,S_{i_0-1}f_2)}{S_{i_0-2}(u+1)g(u+1)}=
\frac{S_{i_0-1}(u)S_{i_0-1}(u+1)}{S_{i_0-2}(u+1)}\;
\frac{\Wrd(f_1,f_2)}{g(u+1)} = R(u) \frac{\Wrd(f_1,f_2)}{g(u+1)}.
\eea
Since the space spanned by $w_1,\dots,w_{i_0-2}\in U$ is acceptable,
the functions $g(u+1)=\Wrd(w_1,\dots,w_{i_0-2})(u+1)/S_{i_0-2}(u+1)$
and $R(u)$ do not have common zeros. Therefore, the discrete Wronskian
$\Wrd(w_1,\dots,w_{i_0})$ is divisible by $R(u)$.
\end{proof}

\subsection{Discrete exponents}
For $\la\in\C$, there exists an increasing sequence of
nonnegative integers $\{c_1 < \dots < c_M\}$
and a basis $\{f_1,\dots,f_M\}$ of $U$
such that for $i=1, \dots, M$, we have
$f_i(\la+j) = 0$ for $j=0, \dots, c_i - 1$ and
$f_i(\la + c_i )\neq 0$.
This sequence of integers is defined uniquely and will be called
{\it the sequence of discrete exponents of $U$ at $\la$}. We say that the basis
$\{f_1,\dots,f_M\}$ agrees with exponents at $\la$.

We say that $\la$ is {\it a singular point of} $U$ if the discrete exponents at
$\la$ differ from the sequence $ \{0< 1 <\dots < M-1\}$.

For $i=1,\dots,M$, introduce {\it the local frame-type polynomials}
\bea
Q_i\ =\ \prod_{j=c_{i-1}-i+2}^{c_i-i}\!\!(u-\la - j)
\eea
where $c_0=-1$. Notice that
\begin{enumerate}
\item[$\bullet$]
roots of each $Q_i$ are simple,
\item[$\bullet$]
$\deg\,Q_i = c_i-c_{i-1}-1$,
\item[$\bullet$]
sets of roots of different polynomials do not intersect,
\item[$\bullet$]
the union of roots of all $Q_i$ is the sequence $\la, \la+1,\dots,
\la+c_M-M$.
\end{enumerate}

\begin{theorem}
\label{thm discrete exponents}
Let $c_1 \leq \dots\leq c_M$ be the discrete exponents of $U$ at
$\la$. Then
\begin{enumerate}
\item[(i)]
The discrete Wronskian of $U$ is divisible by
$$
\prod_{s=1}^M \prod_{j=s-M}^{c_s-M} (u-\la-j)\ = \
\prod_{i=1}^M \prod_{j=1}^{M-i+1} Q_i(u+j-1)\ .
$$
In particular, the total degree of this divisor is
$\sum_{i=1}^M(c_i-i+1)$.
\item[(ii)]
If $S_1, \dots, S_M$ are the polynomials from Lemma \ref{U=T}, then for any $k$,
the polynomial $S_k$ is divisible by
\bea
\prod_{i=0}^{k}\prod_{j=1}^{k-i+1}\, Q_i(u+j-1) \ .
\eea
\hfill $\square$
\end{enumerate}
\end{theorem}

\begin{cor}
\label{cor divisibility in Q}
Assume that the exponents of $U$ at $\la$ have the form
\bea
0 < \dots < M-L-1 < M-L + l_{1}< M-L+l_{2}< \dots < M-L+l_{L}
\eea
for a suitable $L$, $1\leq L\leq M$, and $0<l_1 < \dots < l_L$. Then
the local frame-type polynomials have the form
$Q_i=1$ for $i=1,\dots, M-L$, \
\bea
Q_{M-L+1} = \prod_{j=0}^{l_1-1}(u-\la - j)\ ,
\qquad
Q_{M-L+i} = \prod_{j=l_{i-1}-i+2}^{l_i-i}(u-\la - j)
\eea
for
$i=2,\dots,L$, and the discrete Wronskian of $U$ is
divisible by
\bea
\prod_{k=1}^L \prod_{j=k-L}^{l_k-L} (u-\la-j)\ .
\eea
In particular, the total degree of this divisor is $\sum_{k=1}^L(l_k-k+1)$.
\end{cor}

\subsection{Proof of Theorem \ref{thm discrete exponents}}
We shall prove part(i). Part (ii) is proved analogously.

Let $\{f_1, \dots ,f_M\}$ be a basis in $U$ that agrees with exponents
at $\la$. Consider the matrix-valued function $F(u) = [ f_j(u+k-1)
]_{j,k=1,\dots,M}$\ .

\begin{lem}
\label{lem 1 corank}
For $t\in \C$, if the corank of $F(t)$ is $r$,
then the discrete Wronskian
\linebreak
$\Wrd(f_1,\dots,f_M)$ is divisible by $(u-t)^r$. 
\hfill $\square$
\end{lem}

For $j\in\Z$, set $d_{j} = \#\{ s \le M\ |\ j+M \leq c_s \}$.
For $j\geq 0$, set $r_{j}=d_j$, and for $j<0$,
set $r_{j}=\max\,(0,d_j+j)$.

It is easy to see that for $j< 0$, the number $r_{j}$
can be also defined as $\#\{ s \le M\ |\ s\leq j+M\leq c_s \}$.

\begin{lem}
\label{lem 2 corank}
For $j\in\Z$, the corank of $F(\la+j)$ is not less than $r_{j}$.
\end{lem}
\begin{proof}
If $j \ge 0$, then $F(\la+j)$ has $d_{j}$ zero rows
produced by the functions $f_{M-d_{j}+1}, \dots ,f_M$.
Hence, the corank of $F(\la+j)$ is at least $ d_{j} = r_{j}$.

If $-M\leq j<0$, then the rows produced by the functions
$f_{M-d_{j}+1}, \dots ,f_M$ have zeros everywhere except in the first $-j$ columns.
Hence, the corank of $F(\la+j)$ is at least $d_{j} - (-j)$.
\end{proof}

By Lemma \ref{lem 2 corank}, the discrete Wronskian of $U$ is
divisible by $\prod_{j=1-M}^{c_M-M} (u-\la-j)^{r_{j}}$. That can be
written as $\prod_{s=1}^M \prod_{j=s-M}^{c_s-M}
(u-\la-j)$. Theorem \ref{thm discrete exponents} is proved.

\subsection{Numerically non-degenerate space of quasi-exponentials}
\label{section non-degenerate space}
Let $U$ be a space of quasi-expo\-nen\-tials like in
Section \ref{subsec quasi-exponentials}.

Let $\la_1,\dots,\la_n \in \C$ be such that $\la_i-\la_j \notin \Z$
for $i\neq j$. For $i=1,\dots, n$, let the exponents of $U$ at
$\la_i$ have the form
\bea
\{ 0 < \dots < M-N_i-1 < M-N_i + n_{i1}<
M-N_i+n_{i2}< \dots < M-N_i+n_{iN_i} \}
\eea
for a suitable $N_i$, $1\leq N_i\leq M$, and $0<n_{i1}<\dots < n_{iN_i}$.
Set
\bea
N\ =\ N_1\ +\ \dots\ +\ N_n\ .
\eea

We say that the space $U$ is {\it a space of quasi-exponentials with data}\
\bea
\frak D \ =\ \{\, m, M_a, m_{ab}, z_a,\, n, N_i, n_{ij}, \la_i\, \}
\eea
where $a=1,\dots,m,\, b=1,\dots,M_a$,\, $i=1,\dots, n,\, j=1,\dots,N_i$.

We say that the space $U$ is {\it a numerically
non-degenerate space of quasi-exponentials with respect to data} $\frak D$\ if
\bean
\label{numerically non-degenerate}
\sum_{a=1}^m\sum_{b=1}^{M_a}\ (\,m_{a b} +1-b\,)
\ = \ \sum_{i=1}^n\sum_{j=1}^{N_i} \ (\,n_{i j}+1-j\,) \ .
\eean

\medskip

For $i=1,\dots,n$, let $Q_{i1},\dots,Q_{iM}$ be the local frame-type
polynomials associated with the point $\la_i$.
For $k=1,\dots,M$, define
\bean
\label{Q_k}
Q_k(u)\ = \ \prod_{i=1}^n\, Q_{ik}(u)\ .
\eean

\begin{lem}
\label{lem exponents-frame}
If $U$ is a numerically non-degenerate space of quasi-exponentials,
then the discrete Wronskian of $U$ is given by the formula
\bea
\Wrd_U\ = \ \prod_{a=1}^m z_a^{M_au}\
\prod_{k=1}^{M}\prod_{j=1}^{M-k+1}\! Q_k(u+j-1) \ .
\eea
Moreover, if $S_1, \dots, S_M$ are the polynomials from Lemma
\ref{U=T}, then for any $i$,
the polynomial $S_i$ is divisible by
\bea
\prod_{k=0}^{i}\prod_{j=1}^{i-k+1}\, Q_k(u+j-1) \ .
\eea
\end{lem}

The lemma follows from Lemmas \ref{lem on wronskian}, \ref{U=T}, and
Theorem \ref{thm discrete exponents}.
\medskip

\subsection{Fundamental difference operator}
\label{subsection fund difference oper}
{\it The monic fundamental difference operator of a space
of quasi-exponentials $U$} is the unique monic linear
difference operator
\bea
\bar D_U\ =\ \tau_u^M + \bar B_1\tau_u^{M-1} + \dots + \bar B_{M-1}\tau_u + \bar B_M \
\eea
of order $M$ annihilating $U$. Here
\bea
\bar B_i\ =\ (-1)^i\, \frac {\Wrd_{U,\,i}} {\Wrd_U}\ ,
\eea
where $\Wrd_U$ is the discrete Wronskian of a basis $\{f_1,\dots ,f_M\}$
of $ U$, \ {} $\Wrd_{U,\,i}$ is the determinant of the $M\times M$-matrix
whose $j$-th row is
$$
f_j(u), \ f_j(u+1),\ \dots ,\
f_j(u+N-i-1),\ f_j(u+(N-i+1),\
\dots , \ f_j(u+N)\ .
$$

Clearly, $\bar B_1,\dots,\bar B_M$ are rational functions.

\begin{lem}
\label{coeff lemma}
For any $i$, the function $\bar B_i$ has a limit
as $u$ tends to infinity. Denoted this limit by
$\bar B_i(\infty)$. Then
\bea
\phantom{aaaaaa}
x^M\ +\ \bar B_1(\infty)\, x^{M-1}\ +\ \dots\ +\ \bar B_M (\infty)\ =\
\prod_{a=1}^M \, (x-z_a)^{M_a}\ .
\phantom{aaaaaaaaaaaaaaa}
\hfill
\square
\eea
\end{lem}

\begin{lem}
\label{important}
For any $i=1, \dots, M$, the function
\bea
\bar B_i(u)\ \prod_{i=1}^n\prod_{j=1}^{N_i} (u - \la_i-n_{ij}+N_i)
\eea
is a polynomial.
\end{lem}

\begin{proof}
Let $\la$ be one of the points of the set $\{\la_1,\dots,\la_n\}$.
For such a $\la$, in
the proof of Theorem \ref{thm discrete exponents},
we defined the numbers
$d_j$ and $r_j$ for $j\in \Z$.

For $j\in \Z$,
define the new numbers $p_j$ as follows.
Set $p_j=d_{j+1}$ for $j\ge 0$, and set $p_j=\max(0,j+d_{j+1})$
for $j<0$.

By the construction, we have $p_j \le r_j$.

The reasons in the proof of Theorem
\ref{thm discrete exponents} show that
the discrete Wronskian $\Wrd_U$ is divisible by
\bea
X(u) = \prod_{j=1-M}^{c_M-M} (u-\la-j)^{r_j}\ .
\eea
Similar reasons show that for any $k$, the
determinant $\Wrd_{U,k}$ is divisible by
\bea
Y(u) = \prod_{j=1-M}^{c_M-M-1} (u-\la-j)^{p_j}\ .
\eea
As in the end of the proof of Theorem \ref{thm discrete exponents},
we have
\bea
X(u) = \prod_{s=1}^M \prod_{j=s-M}^{c_s-M} (u-\la-j) =
\prod_{\fratop{s=1}{c_s\ge s}}^M \ \prod_{j=s-M}^{c_s-M} (u-\la-j)\ ,
\eea
where in the second expression we excluded the empty products over $j$.
Similarly,
\bea
Y(u) = \prod_{\fratop{s=1}{c_s\ge s+1}}^M \prod_{j=s-M}^{c_s-1-M} (u-\la-j) =
\prod_{\fratop{s=1}{c_s\ge s}}^M \prod_{j=s-M}^{c_s-1-M} (u-\la-j)\ ,
\eea
where the first expression is analogous to the second expression for $X(u)$,
and the second expression may contain certain empty products over $j$.

Using the second expressions for $X(u)$ and $Y(u)$, we get
\bean
\label{X=Y}
X(u)\ =\ Y(u)\, \prod_{\fratop{s=1}{c_s\ge s}}^M \,(u-\la-c_s+M)\ .
\eean

Now if $\la=\lambda_i$, we shall provide $X(u)$ and $Y(u)$ with index $i$.
Calculating the product in \Ref{X=Y} for $\la=\la_i$,
we get
\bea
X_i(u)\ =\ Y_i(u)\, \prod_{j=1}^{N_i}\, (u-\la_i-n_{ij}+N_i)\ .
\eea
Multiplying this formula over $i=1,\dots,n$, we conclude that for any $k$,
the product
\bea
\Wrd_{U,k}(u)\ \prod_{i=1}^n \prod_{j=1}^{N_i}\, (u-\la_i-n_{ij}+N_i)
\eea
is divisible by the discrete Wronskian $\Wrd_U$.
\end{proof}

Define {\it the regularized fundamental difference operator
of the space} $U$ as the linear difference operator
\bea
D_U\ =\ B_{0}(u) \,\tau_u^M\, +\, B_{1}(u)\, \tau_u^{M-1}\, +\,
\dots\, +\, B_{M}(u)\
\eea
of order $M$ with
polynomial coefficients, which annihilates $U$ and such that its leading
coefficient $B_0$ is a monic polynomial of the minimal possible
degree.

We have
\bea
\deg\, B_0\ = \ \deg\, B_M\ \ge\ \deg\, B_i\
\qquad
{\rm for}
\quad
i = 1 ,\dots , M-1\ ,
\eea
by Lemma \ref{coeff lemma}.

We say that the space $U$ with data
$\{\, m, M_a, m_{ab}, z_a,\, n, N_i, n_{ij}, \la_i\, \}$
is {\it a
non-degenerate space of quasi-exponentials} if U
is numerically non-degenerate and
\bea
\deg\, B_0\ = \ N\ =\ N_1 + \dots + N_n\ .
\eea

By Lemma \ref{important}, if the space $U$ with data
$\{\, m, M_a, m_{ab}, z_a,\, n, N_i, n_{ij}, \la_i\, \}$ is
non-degenerate, then
\bea
B_0\ =\ \prod_{i=1}^n\prod_{j=1}^{N_i} \ (u - \la_i-n_{ij}+N_i)\ .
\eea

\medskip
\noindent
{\bf Example.}
Let $U$ be the vector space spanned by the quasi-exponentials
$u$ and $u(u-1)$. Then $M=2, \,m=1,\, z_1=1,\, m_{11}=1,\, m_{12}=2$.
Let $\la_1=0$.
The exponents of $U$ at $\la_1$ are $1$ and $ 2$. Then
$n_{11}=1,\, n_{12}=2,\, N = N_1 = 2$.
Equality \Ref{numerically non-degenerate}
takes the form: 2=2. We have
\bea
D_U\ =\
u(u+1)\,\tau_u^2\, - \,2 u(u+2)\,\tau_u \,+\, (u+1)(u+2) \ .
\eea
Hence $U$ with this data is non-degenerate.

\medskip

\noindent
{\bf Example.}
Let $U$ be the vector space spanned by the quasi-exponentials
$u$ and $(-1)^uu$. Then $M=2,\, m=2,\, z_1=0,\, z_2=-1,
\,M_1=1,\, m_{11}=1,\, M_2=1,\, m_{21}=1$.
Let $\la_1=0$.
The exponents of $U$ at $\la_1$ are $1$ and $ 2$. Then
$n_{11}=1,\, n_{12}=2, N=N_1=2$.
Equality \Ref{numerically non-degenerate}
takes the form: 2=2. Hence $U$ with this data is numerically non-degenerate.
We have
\bea
D_U\ =\
u\,\tau_u^2\,+\, (u+2) \ .
\eea
Hence $U$ with this data is degenerate.

\medskip

\noindent
{\bf Example.}
Let $U$ be the vector space spanned by the quasi-exponentials
$u$ and $(-1)^uu$. Then $M=2,\, m=2,\, z_1=0,\, z_2=-1,
\,M_1=1,\, m_{11}=1,\, M_2=1,\, m_{21}=1$.
Let $\lambda_1=-1$. The exponents of $U$ at $\lambda_1$ are $0,3$. Then
$N_1=1$,\, $n_{11}=2$,\, $N=1$ and Equality \Ref{numerically non-degenerate}
is $2=2$. Therefore, $U$ with this
data is numerically non-degenerate. With this data, we have
\bea
D_U \ =\ u \,\tau_u^2\ +\ (u+2)\ .
\eea
Hence, $U$ with this data is non-degenerate.

\begin{theorem}
\label{theorem on reg difference operator}
Assume that the space $U$ with data
$\{\, m, M_a, m_{ab}, z_a,\, n, N_i, n_{ij}, \la_i\, \}$
is a non-degenerate space of quasi-exponentials.
Then
\begin{enumerate}
\item[(i)]
We have
\bea
B_M\ =\ \prod_{a=1}^m (-z_a) ^{M_a}
\prod_{i=1}^n\prod_{j=1}^{N_i} \ (u - \la_i + j)\ .
\eea
\item[(ii)]
Write 
$$
D_U\ =\ u^N A_0(\tau_u) + u^{N-1} A_1(\tau_u) + \dots + A_N(\tau_u)
$$
where $A_i(\tau_u)$ is a polynomial in $\tau_u$ with constant coefficients. Then
\bea
&
A_0(\tau_u) \ =\ \prod_{a=1}^m (\tau_u - z_a)^{M_a}\ .
\eea
\item[(iii)]
The polynomials $A_0,\dots,A_M$ have no common factors of positive degree.
\end{enumerate}
\end{theorem}

\begin{cor}
If the space $U$ is non-degenerate with respect to a data
$\{\, m, M_a, m_{ab},$ $z_a,\, n, N_i, n_{ij}, \la_i\, \}$, then the data
is determined uniquely.
\end{cor}

The corollary follows from part (i) of the theorem.

\medskip
\noindent
{\it Proof of Theorem \ref{theorem on reg difference operator}.}\
Part (ii) follows from Lemma \ref{coeff lemma}.
Part (iii) follows from the fact that $U$ does not contain exponential functions
$z^u$.

Let $Q_1,\dots,Q_M$ be the polynomials introduced in \Ref{Q_k}.
To prove part (i) it is enough to notice that
\bea
\frac{B_M(u)}{B_0(u)}
&=& (-1)^M
\frac {\Wrd (u+1)}{\Wrd(u)} = \prod_{k=1}^M \frac {Q_k(u+M+1-k)} {Q_k(u)}
\ \prod_{a=1}^m z_a^{M_a}
\\
&=& (-1)^M
\prod_{i=1}^n\prod_{j=1}^{N_i} \ \frac
{u - \la_i + j}
{u - \la_i-n_{ij}+N_i}
\ \prod_{a=1}^m z_a^{M_a} \ .
\eea
\hfill
$\square$

\subsection{Regularized Conjugate Space}
\label{Subsection conj space}
Let $U$ with data
$\{\, m, M_a, m_{ab}, z_a,\, n, N_i, n_{ij}, \la_i\, \}$
be a non-degenerate space of quasi-exponentials as in Section
\ref{subsection fund difference oper}. Let $\Wrd_U$ be the discrete
Wronskian of $U$ and let $B_M(u)$ be the last coefficient of the
regularized fundamental difference operator of $U$.

The complex vector space spanned by all functions of the form
\bean
\label{con space}
\frac
{\tau_u( \Wrd(f_1,\dots,f_{M-1}))}
{B_M(u)\,\Wrd_U(u)}
\eean
with $f_i\in V$ has dimension $M$. This space is denoted by
$U^\ddagger$
and called {\it regularized conjugate to} $U$.

\begin{lem}
\label{lem on poles}
For any $g \in U^\ddagger$, the function
\bea
g(u)\ \prod_{i=1}^n\prod_{j=0}^{n_{iN_i}} (u-\la_i + N_i- j)
\eea
is holomorphic in $\C$.
\end{lem}
\begin{proof}
Let $Q_1,\dots,Q_M$ be the polynomials introduced in \Ref{Q_k}.
Let $g$ be a function in \Ref{con space}.
Then ${\tau_u(\Wrd(f_1,\dots,f_{M-1}))}$ is divisible by
$$
\prod _{i=1}^{M-1}Q_1(u+i)\cdot \prod
_{i=1}^{M-2}Q_2(u+i) \cdots\ Q_{M-1}(u+1)
$$
while
$$
\Wrd_U(u) =
\prod_{a=1}^m z_a ^{M_a} \prod _{j=1}^{M} Q_1(u+j-1)\cdot \prod
_{j=1}^{M-1}Q_2(u+j-1)\ \cdots\ Q_{M}(u)\ .
$$
Hence, the possible poles of
$g$ come from the product $B_M(u) Q_1(u)\dots Q_M(u)$ which remains
in the denominator of $g$. But this product is exactly the product in
Lemma \ref{lem on poles}.
\end{proof}

For every $i =1, \dots, n$ and $j = 1, \dots, N_i$,
fix a function $g_{ij}$ in $U$
which is equal to zero at $u=\la_i, \la_i+1, \dots, \la_i + M-N_i + n_{ij}-1$
and which is not equal to zero at $u=\la_i + M-N_i + n_{ij}$.

\begin{lem}
\label{special lem on poles}
For given $i = 1, \dots, n$, \, $j = 0, \dots, N_i-1$, let
$f_1,\dots,f_{M-1}$ be a collection of functions in $U$ containing the functions
$g_{i,j+1},g_{i,j+2}, \dots, g_{i,N_i}$ and let
\bea
F(u)\ =\
\frac
{\tau_u( \Wrd(f_1,\dots,f_{M-1}))}
{B_M(u)\,\Wrd_U(u)} \ .
\eea
Then for $j=0$, the function $F$ has no poles at
$$
u\ =\ \la_i-N_i,\ \la_i-N_i+1,\ \dots,\ \la_i-N_i+n_{iN_i}\ .
$$
For $j>0$, the function
$F$
has no poles at
$$
u\ =\ \la_i-N_i+n_{ij}+1, \ \la_i-N_i+n_{ij}+2,\ \dots,\ \la_i-N_i+n_{iN_i}\ .
$$
Moreover, if
$f_1,\dots,f_{M-1}$ is a generic collection of functions in $U$ containing the functions
$g_{i,j+1},g_{i,j+2}, \dots, g_{i,N_i}$, then the function $F$
has a nonzero residue at $u=\la_i-N_i+n_{ij}$.
\end{lem}

\begin{proof}
The first two statements of the lemma are proved in the same way as
Lemma \ref{lem on poles}.

We shall prove that the residue of $F$ at $u=\la_i-N_i+n_{ij}$ is nonzero first
assuming that $n_{ij} \geq N_i$. To prove that the residue is nonzero it is
enough to show that
\beq
\label{ord=1+ord}
\qquad{\rm ord}_{u=\la_i-N_i+n_{ij}}\Wrd_U =
1+{\rm ord}_{u=\la_i-N_i+n_{ij}+1}\Wrd_U = N_i-j+1\ ,
\eeq
and
\beq
\label{ord=ord}
\qquad{\rm ord}_{u=\la_i-N_i+n_{ij}+1}\Wrd_U =
{\rm ord}_{u=\la_i-N_i+n_{ij}+1}
\Wrd (f_1, \dots,f_{M-1})\ .
\eeq
Here ${\rm ord}_{u=\la} f$ denotes the order of zero of the function $f$
at $u=\la$.

Equalities \Ref{ord=1+ord} follow from the numerical non-degeneracy of
the space $U$ and the condition $n_{ij}\ge N_i$.

Since the collection $f_1,\dots ,f_{M-1}$ is a generic collection containing
the functions $g_{i,j+1}$, \dots,  $g_{i,N_i}$, the discrete Wronskian
$\Wrd(f_1,\dots, f_{M-1},g_{ij})$ is nonzero and proportional to $\Wrd_U$.
Expanding the determinant $\Wrd(f_1,\dots, f_{M-1}, g_{ij})$ with respect to
the last row, we have
\begin{align}
\label{ORDER}
\qquad\Wrd_U (u)\ =\ {\rm const}\ \bigl(\,
g_{ij}(u+M-1) \ \Wrd (f_1, \dots, f_{M-1})(u)\ - \ \dots &
\\
{}-\;(-1)^M g_{ij}(u)\ \tau_u (\Wrd (f_1, \dots, f_{M-1}))(u)\,\bigr)\! & \,\ .
\notag
\end{align}
where ${\rm const}\ne 0$.

The order of $\Wrd (f_1, \dots,f_{M-1})$ at $u=\la_i-N_i+n_{ij}+1$ is not less
than $N_i-j$. This follows from Theorem \ref{thm discrete exponents} applied to
the functions $f_1,\dots, f_{M-1}$. A similar reason shows that the order at
$u=\la_i-N_i+n_{ij}+1$ of all of the other $(M-1)\times (M-1)$ minors in
the right hand side of \Ref{ORDER} is also not less than $N_i-j$.

By the construction, the function $g_{ij}$ is nonzero at $u=\la_i-N_i+n_{ij}$
and is zero at $u=\la_i-N_i+n_{ij}-l$ for $l=1,\dots,M-1$. Therefore, the only
term in the right hand side of equality \Ref{ORDER} that can have order $N_i-j$
at $u=\la_i-N_i+n_{ij}+1$ is $g_{ij}(u+M-1)\;\Wrd(f_1,\dots,f_{M-1})(u)$.
Since ${\rm ord}_{u=\la_i-N_i+n_{ij}+1}\Wrd_U = N_i-j$, this shows that
the orders of $\Wrd_U$ and $\Wrd (f_1, \dots,f_{M-1})$ at $u=\la_i-N_i+n_{ij}$
are equal.

\smallskip
To prove that the residue of $F$ at $u=\la_i-N_i+n_{ij}$ is nonzero in the case
$n_{ij} < N_i$, it is enough to show that
$$
\qquad{\rm ord}_{u=\la_i-N_i+n_{ij}}\Wrd_U =
{\rm ord}_{u=\la_i-N_i+n_{ij}+1}\Wrd_U = n_{ij}-j+1\ ,
$$
and
$$
\qquad{\rm ord}_{u=\la_i-N_i+n_{ij}+1}\Wrd_U =
{\rm ord}_{u=\la_i-N_i+n_{ij}+1}
\Wrd (f_1, \dots,f_{M-1})\ .
$$
The proof is similar to the proof of equalities \Ref{ord=1+ord} and
\Ref{ord=ord}.
\end{proof}

\begin{theorem}
\label{Theorem on difference operator fo dagger}
If $D_U = B_{0}(u) \tau_u^M + B_{1}(u) \tau_u^{M-1} +
\dots + \tau_u B_{M-1}(u) + B_{M}(u)$
is the regularized fundamental difference operator of $U$, then the operator
\bea
D_U^\ddagger \ =\
\tau_u^M B_{M}(u) + \tau_u^{M-1} B_{M-1}(u) +
\dots + \tau_u^{1} B_{1}(u) + B_{0}(u)
\eea
annihilates $U^\ddagger$.

\end{theorem}

\begin{proof} Consider the scalar
equation $\frac 1{B_0(u)}D_U y = 0$ with respect to an unknown
function $y(u)$. For $i=1,\dots,M$, introduce $w_i = \tau_u^{i-1}y$
and present the equation as a system of first order equations
\bea
\tau_u
w_{M} = -\frac {B_1(u)}{B_0(u)} w_{M}- \dots - \frac {B_M(u)}{B_0(u)}
w_{1}\ ,
\qquad
\tau_u w_i = w_{i+1}
\eea
for $i = 1,\dots, M-1$. For the column
vector $w = (w_1,\dots,w_M)$, the system can be presented as a matrix
equation $\tau_u w = C w$ with the $M\times M$-matrix
\[
C =
\left( \begin{array} {cccccc}
0 & 1 & 0 & \dots\dots& 0 & 0
\\
0 & 0 &{} 1 & \dots\dots & 0 & 0
\\
{}\dots & {} \dots & \dots & \dots\dots& {} \dots &{} \dots
\\
{}0 & 0 & 0 & \dots\dots & 0 & 1
\\
- \frac {B_M(u)}{B_0(u)} & - \frac {B_{M-1}(u)}{B_0(u)}& - \frac {B_{M-2}(u)}{B_0(u)}
& \dots\dots &
- \frac {B_2(u)}{B_0(u)} & - \frac {B_1(u)}{B_0(u)}
\end{array} \right).
\]
Let $\Psi$ be a fundamental $M\times M$-matrix of solutions,\
$\tau_u \Psi = C \Psi$. Then $\tau_u \Psi^{-1} = \Psi^{-1} C^{-1}$ where
\[
C^{-1} =
\left( \begin{array} {cccccc}
- \frac {B_{M-1}(u)}{B_M(u)} & - \frac {B_{M-2}(u)}{B_M(u)}& - \frac {B_{M-3}(u)}{B_M(u)}
& \dots\dots &
- \frac {B_1(u)}{B_M(u)} & - \frac {B_0(u)}{B_M(u)}
\\
1 & 0 & 0 & \dots\dots& 0 & 0
\\
0 & 1 & 0 & \dots\dots & 0 & 0
\\
\dots & \dots & \dots & \dots\dots& \dots & \dots
\\
0 & 0 & 0 & \dots\dots & 1 & 0
\end{array} \right).
\]
For a row vector $(v_1,\dots,v_M)$ the equation $\tau_u v = v C^{-1}$
has the form:
\bea
\tau_u v_M = - \frac {B_{0}}{B_{M}} v_1,
& \qquad &
\tau_u v_{M-1} = - \frac {B_{1}}{B_{M}} v_1 + v_M,
\\
\tau_u v_{M-2} = - \frac {B_{2}}{B_{M}} v_1 + v_{M-1}, \quad & \dots\ ,&
\quad
\tau_u v_{1} = - \frac {B_{M-1}}{B_{M}} v_1 + v_{2}\ .
\eea
This system reduces to the scalar equation $D_U^\star v_1 = 0$ where
\bea
D^\star_U = \tau_u^{M} + \tau_u^{M-1}\frac {B_{M-1}(u)}{B_M(u)} + \dots
+ \tau_u \frac {B_1(u)}{B_M(u)} + \frac {B_0(u)}{B_M(u)} \ .
\eea
Thus the kernel of the difference operator $D^\star_U$
consists of the first row entries of the
matrix $\Psi^{-1}$.

If $\{f_1,\dots,f_M\}$ is a basis of $U$, then
\[
\Psi =
\left( \begin{array} {cccc}
f_1(u) & f_2(u) & \dots\dots& f_M(u)
\\
f_1(u+1) & f_2(u+1) & \dots\dots& f_M(u+1)
\\
\dots & \dots & \dots\dots& \dots
\\
f_1(u+M-1) & f_2(u+M-1) & \dots\dots& f_M(u+M-1)
\end{array} \right) \ .
\]
The formula for the inverse matrix elements shows that $D^\star_U$
annihilates the functions of the form ${\tau_u(
\Wrd(f_1,\dots,f_{M-1}))}/ {\Wrd_U(u)}$. Then the operator
$D^\ddagger_U$ annihilates the functions of the form \Ref{con space},
since $D^\ddagger_U = D^\star_U \cdot B_M$ where \, $\cdot B_M$ is the
operator of multiplication by the function $B_M$.
\end{proof}

\section{Integral Transforms}
\label{Sec Int transforms}

\subsection{Mellin-type transform}\label{sec int transform}
Let $V$ be a non-degenerate space of quasi-polynomials
with data
$\frak D_V = \{n, N_i,n_{ij},\la_i,$ $m, M_a, m_{ab},$ $ z_a\}$
where $i=1,\dots, n,\, j=1,\dots,N_i,\, a=1,\dots,m,\, b=1,\dots,M_a$.
Let $V^\dagger$ be the space regularized conjugate to $V$.

For $a = 1, \dots , m$, denote by $\gamma_a$ a small
circle around $z_a$ in $\C$ oriented counterclockwise.

Denote by $U$
the complex vector space spanned by functions of the form
\bean
\label{int transform}
\hat f_a(u)\ =\ \int_{\gamma_a}\, x^{u}\,
f(x)\, dx\ ,
\eean
where $a = 1, \dots , m$, \ $f \in V^\dagger$. The
vector space $U$ is called {\it bispectral dual to $V$}.

\begin{theorem}
\label{thm Mellin transform}
Let $V$ be a non-degenerate space of quasi-polynomials
with data
$$
\frak D_V \ =\ \{\,n, N_i,n_{ij},\la_i,\, m, M_a, m_{ab}, z_a\,\}
$$
where $i=1,\dots, n,\, j=1,\dots,N_i,\, a=1,\dots,m,\, b=1,\dots,M_a$. Then

\begin{enumerate}
\item[(i)] The space $U$ is a non-degenerate space of quasi-exponentials with data
$$
\frak D_U \ =\ \{\, m, M_a, m_{ab}, z_a, \,n, N_i, n_{ij}, \la_i + N_i\, \}
$$
where $a=1,\dots,m,\, b=1,\dots,M_a$,\, $i=1,\dots, n,\, j=1,\dots,N_i$.

\item[(ii)]
Let $D_V = \sum_{i=1}^M\sum_{j=1}^N
A_{ij}x^i (x\partial_x)^j$ be the regularized fundamental differential
operator of $V$
where $A_{i j}$ are suitable complex numbers. Then
\bea
\sum_{i=1}^M\sum_{j=1}^N\, A_{i j}\ u^j \,\tau_u^i
\eea
is the regularized fundamental difference operator of $U$.
\end{enumerate}
\end{theorem}

The theorem is proved in Section \ref{PTMT}.

\subsection{Fourier-type transform}
\label{sec F transform}
Let $U$ be a non-degenerate space of quasi-exponentials
with data $\frak D_U = \{ m, M_a, m_{ab}, z_a, n, N_i, n_{ij}, \la_i \}$
where $a=1,\dots,m,\, b=1,\dots,M_a$,\, $i=1,\dots, n,\, j=1,\dots,N_i$.
Let $U^\ddagger$ be the space regularized conjugate to $U$.

For $i = 1,\dots , n$, consider the arithmetic sequence
$$
\frak S_i\ =\ \{\,\la_i-N_i,\la_i-N_i+1,\la_i-N_i+2, \dots , \la_i-N_i+n_{iN_i}\,\}
$$
consisting of $n_{iN_i}+1$ terms.

For $i=1,\dots, n$, fix a non-self-intersecting closed connected curve
$\gamma_i$ in $\C$ oriented counterclockwise and such that the points of
the sequence $\frak S_i$ are inside $\gamma_i$ and the points of other
sequences $\frak S_j$ for $j\neq i$ are outside $\gamma_i$.

Denote by $V$
the complex vector space spanned by functions of the form
\bean
\label{int F transform}
\hat f_i(x)\ =\ \int_{\gamma_i}\, x^{u}\,
f(u)\, du \ ,
\eean
where $i = 1, \dots , n$, \ $f \in U^\ddagger$. The
vector space $V$ is called {\it bispectral dual to $U$}.

\begin{theorem}
\label{thm F transform}
Let $U$ be a non-degenerate space of quasi-exponentials
with data
$$
\frak D_U \ =\ \{\, m, M_a, m_{ab}, z_a,\, n, N_i, n_{ij}, \la_i\, \}
$$
where $a=1,\dots,m,\, b=1,\dots,M_a$,\, $i=1,\dots, n,\, j=1,\dots,N_i$. Then

\begin{enumerate}
\item[(i)] The space $V$ is a non-degenerate space of quasi-polynomials with data
$$
\frak D_V\ =\ \{\,n, N_i,n_{ij},\la_i-N_i,\, m, M_a, m_{ab}, z_a\,\}
$$
where $i=1,\dots, n,\, j=1,\dots,N_i,\, a=1,\dots,m,\, b=1,\dots,M_a$.

\item[(ii)]
Let $D_U = \sum_{i=1}^N\sum_{j=1}^M
A_{ij}u^i \tau_u^j$ be the regularized fundamental difference
operator of $U$
where $A_{i j}$ are suitable complex numbers. Then
\bea
\sum_{i=1}^N\sum_{j=1}^M \,A_{i j}\ x^j\, (x\partial_x)^i
\eea
is the regularized fundamental differential operator of $V$.
\end{enumerate}
\end{theorem}

The theorem is proved in Section \ref{PTFT}.

\medskip

Theorems \ref{thm Mellin transform}
and \ref{thm F transform}
imply that if $V$ is a non-degenerate space of quasi-polynomials and $U$ is bispectral dual
to $V$, then $V$ is bispectral dual to $U$. Similarly,
if $U$ is a non-degenerate space of quasi-exponentials and $V$ is bispectral dual
to $U$, then $U$ is bispectral dual to $U$.

\subsection{Proof of Theorem \ref{thm Mellin transform}}
\label{PTMT}
The exponents of $V^\dagger$ at $z_a$ are
$$
\{-m_{a,M_a}-1 < \dots
< -m_{a1}-1< 0< 1<\dots < N - M_a -1 \}\ .
$$
Integral \Ref{int transform} is nonzero only if $f$ has a pole at
$x=z_a$. If $f$ has a pole at $x=z_a$ of order $-m_{ab}-1$, then
integral \Ref{int transform} has the form $z_a^u q_{ab}(u) $ where
$q_{ab}$ is a polynomial in $u$ of degree $m_{ab}$. Thus $U$ is a
space of quasi-exponentials of dimension $M$ generated by
quasi-exponentials $z_a^u q_{ab}(u) $ where $a=1,\dots,m,$\,
$b=1,\dots,m_{ab}$ and $q_{ab}$ is a polynomial in $u$ of degree
$m_{ab}$.

It is clear that the operator
$D^\dagger = \sum_{i=1}^M\sum_{j=1}^N\ A_{i j}\ u^j \tau_u^i$\
annihilates $U$. Write
\bea
D^\dagger\ =\ B_0(u)\tau_u^M + \dots + B_{M-1}(u) \tau_u + B_M(u)
\eea
where $B_a(u)$ are polynomials in $u$ with constant coefficients.
Lemma \ref{deg coeff} implies that
$B_0(u)= \prod_{j=1}^{n}(u-\la_j)^{N_j}$ and
the polynomials $B_0, \dots , B_M$ have no common factor of positive degree.
Therefore, $D^\dagger$ is the
regularized fundamental operator of $U$.

The functions $x^{\la_i}p_{ij}(x)$, $i=1,\dots,n,\
j=1, \dots , N_i$,\ form a basis of $V$. For given $i, j$, order
all basis functions of $V$ except the function
$x^{\la_i}p_{ij}(x)$. Denote by $\Wr_{ij}$ the Wronskian
of this ordered set of $N-1$ functions. The functions
$$
f_{ij}(x)\ =\ \frac{\Wr_{ij}(x)}{ \Wr_V(x)\,x^N\,\prod_{a=1}^m\,(x-z_a)^{M_a}}
$$
form a basis in $V^\dagger$.
Such a function $f_{ij}$ has the form
$x^{-\la_i}r_{ij}$ where $r_{ij}$ is a rational function in $x$.
We have
\bea
{\rm ord}_{x=0}\ r_{ij} = -N_i\ ,
\qquad
{\rm ord}_{x=\infty}\ r_{ij} = -n_{ij}-M-1 \ .
\eea
Consider the following element of $U$,
\bea
F_{ij}(u)\
=\ \sum_{a=1}^m \widehat{(f_{ij})}_a (u)\ =\
\int_{\cup_{a=1}^m \gamma_a}x^u f_{ij}(x)\ dx \ .
\eea
If $u = \la_i + m$ where$ \,N_i\leq m \leq n_{ij}+ M-1$,
then the integrand is a rational function with zero residues at $0$ and $\infty$.
Hence $F_{ij}$ is zero at $u=\la_i +m$ for $m$ from
this arithmetic sequence.
This remark together with Theorem
\ref{thm discrete exponents} proves that $U$ is a non-degenerate space of quasi-exponentials
with data\
$\{ \, m, M_a, m_{ab}, z_a,\, n, N_i, n_{ij}, \la_i + N_i\, \}$.
\hfill
$\square$

\subsection{Proof of Theorem \ref{thm F transform}}
\label{PTFT}
By Lemmas \ref{lem on poles} and \ref{special lem on poles},
for any
$i=1,\dots,n$, the functions $\hat f_i(x)$ in \Ref{int F transform}
have the form $x^{\la_i - N_i}p_{ij}(x)$
where $j=1,\dots,N_i$ and $p_{ij}$ is a polynomial
of degree not greater than $n_{ij}$. Moreover, there exists a function
$x^{\la_i - N_i}p_{ij}(x)$ with $p_{ij}$ of degree exactly equal to $n_{ij}$.

The functions $z_a^{u}q_{ab}(u)$, $a=1,\dots,m,\
b=1, \dots , M_a$,\ form a basis of $U$. For given $a, b$, order
all basis functions of $U$ except the function
$z_a^{u}q_{ab}(u)$ and denote them $f_1,\dots, f_{M-1}$.
The corresponding function
\bea
f_{ab}(u)\ =\ \frac
{\tau_u( \Wrd(f_1,\dots,f_{M-1}))}
{B_M(u)\,\Wrd_U(u)} \ \in \ U^\ddagger
\eea
has the form
$z_a^{-u}r_{ab}(u)$ where $r_{ab}$ is a rational function in $u$ and
${\rm ord}_{u=\infty}\ r_{ab} = M_a-m_{ab}-N-1$.
Consider the following element of $V$,
\bea
F_{ab}(x)\
=\ \sum_{i=1}^n \widehat{(f_{ab})}_i (x)\ =\
\int_{\cup_{i=1}^n \gamma_i} x^u f_{ij}(u)\ du \ .
\eea
If $x = z_a$, then the integrand is a rational function in $u$
which tends to zero as
$u$ tends to infinity. Denote by ${}^{(i)}$ the $i$-th derivative. Then
$F_{ab}^{(i)}(z_a) = 0$ for $i=0,1,\dots, N-M_a+ m_{ab}-1$,
and $F_{ab}^{(N-M_a+ m_{ab})}(z_a) \neq 0$.

This reason proves part (i) of the theorem.

From Theorem \ref{Theorem on difference operator fo dagger} and
formulas for the Fourier-type integral transform, it follows that the differential
operator
$\sum_{i=1}^N\sum_{j=1}^M A_{i j}\ x^j (x\partial_x)^i$ annihilates $V$.
From part (iii) of Theorem \ref{theorem on reg difference operator}, it follows that
this operator
is the regularized fundamental differential operator of $V$.
\hfill
$\square$

\section{Special spaces}
\label{special spaces}

\subsection{Spaces of quasi-polynomials
of an $(\bs n,\bs \la,\bs m,\bs z)$-type}
Let $N$ be a natural number, $N>1$.
Let $\bs n = (n_1, \dots , n_N)$ be a vector of nonnegative integers.

Let $\bs \la = (\la_1, \dots , \la_N) \in \C^N$.
Assume that $\la_i-\la_j \notin \Z$ for $i\neq j$.

For $i=1,\dots , N$, let $p_i\in \C[x]$ be a polynomial of degree $n_i$ such
that $p_i(0) \neq 0$. Denote by $V$ the complex vector space spanned by
functions \ $x^{\la_i}p_{i}(x)$, $i=1,\dots, N$. The dimension of $V$ is $N$.

Let $\bs z = (z_1,\dots ,z_M)$, $M>1$,\ be a subset in $\C$ containing
all singular points of $V$. Assume that for $a = 1,\dots , M$, the
set of exponents of $V$ at $z_a$ has the form
$$
\{ 0 < 1 < \dots < N-2< N-1+m_{a}\}\ ,
\qquad
m_a \geq 0\ .
$$
We have
\bea
&
\sum_{i=1}^N \ n_i\ =\ \sum_{a=1}^M\ m_a\ .
\eea
We call the pair $(V, \bs z)$ {\it a space of the $(\bs n, \bs \la, \bs m, \bs z)$-type} or
{\it a special space of quasi-polynomials}.

Let $\bar D_V$ be the monic fundamental differential operator of $V$.
The operator
\bea
&
\tilde D_V\ = \ x^N \prod_{a=1}^M\ (x - z_a)\ \bar D_V
\eea
is called {\it the special fundamental differential operator\/}
of the special space $(V,\bs z)$.

Write the special fundamental differential operator in the form
\bea
\tilde D_V = A_0(x) (x\p)^N\!+A_1(x)(x\p)^{N-1}\!+ A_2(x)(x\p)^{N-2}\!+\dots +A_N(x).
\eea
Then by Lemma \ref{deg coeff}, all of the coefficients $A_0,\dots, A_N$ are
polynomials in $x$ of degree not greater than $M$.
If we write
$$
\tilde D_V = x^M B_0(x\p) + x^{M-1}
B_1(x\p) + \dots + x B_{M-1}(x\p) + B_M(x\p),
$$
where $B_i$ is a polynomial in $x\p$ with constant coefficients,
then
\bea
B_0 \ =\ \prod_{i=1}^N\, (x\p - \la_i-n_i)\ ,
\qquad
B_M \ =\ (-1)^M\prod_{a=1}^Mz_a\, \prod_{i=1}^N \,(x\p - \la_i)\ .
\eea

\subsection{Spaces of quasi-exponentials of an $(\bs m, \bs z, \bs n, \bs \la)$-type}
Let $M$ be a natural number, $M>1$.
Let $\bs m = (m_1, \dots , m_M)$ be a vector of nonnegative integers.

Let $\bs z = (z_1, \dots , z_M)$ be distinct nozero complex numbers with fixed
arguments.

For $a=1, \dots , M$, let $q_a\in \C[u]$ be a polynomial of degree $m_a$.
Denote by $U$ the complex vector space spanned by the functions \ $z_a^{u}q_a(u)$,\
$a=1,\dots , M$. The dimension of $U$ is $M$.

Assume that for some natural number $N>1$, there exists a
subset of distinct numbers $\bs \la = (\la_1,\dots ,\la_N)$ in $\C$
with three properties:
\begin{enumerate}
\item[$\bullet$]
$\la_i-\la_j \notin \Z$ for $i\neq j$.
\item[$\bullet$]
For $i = 1,\dots , N$, the set of discrete exponents of $U$ at
$\la_i$ has the form
$$
\{ 0 < 1 < \dots < M-2< M-1+n_{i}\}\ ,
\qquad
n_i \geq 0\ .
$$
\item[$\bullet$]\ {}
$\sum_{i=1}^N \, n_i\ =\ \sum_{a=1}^M\, m_a $.
\end{enumerate}
In this case we
call the pair $(U, \bs \la)$ {\it a space of quasi-exponentials
of the $(\bs m, \bs z, \bs n, \bs \la)$-type} or
{\it a special space of quasi-exponentials}.

Let $\bar D_U$ be the monic fundamental difference operator of $U$.
The operator
\bea
&
\tilde D_U\ = \ \prod_{i=1}^N\ (u - \la_i -n_i+1)\ \bar D_V
\eea
is called {\it the special fundamental difference operator } of the special space $(U,\bs \la)$.

Write the special fundamental difference operator in the form
$$
\tilde D_U = B_0(u) \tau_u^M\!+B_1(u)\tau_u^{M-1}\!+ B_2(u)\tau_u^{M-2}\!
+\dots +B_M(u).
$$
Then by Lemma \ref{important} and Theorem
\ref{theorem on reg difference operator},
all of the coefficients $B_0,\dots, B_M$ are
polynomials in $u$ of degree not greater than $N$ and
$$
B_M \ =\ (-1)^M\,\prod_{a=1}^M z_a\, \prod_{i=1}^N \,(u - \la_i+1)\ .
$$
If we write
$$
\tilde D_U = u^N A_0(\tau_u) + u^{N-1}A_1(\tau_u) + \dots+ u A_{N-1}(\tau_u) + A_N(\tau_u),
$$
where $A_i(\tau_u)$ is a polynomial in $\tau_u$ with constant coefficients,
then
\bea
&
A_0(\tau_u) \ =\ \prod_{a=1}^M\, (\tau_u - z_a)\ .
\eea

\subsection{Special Mellin-type transform}
\label{Special Mellin-type transform}
Let $(V,\bs z)$ be a space of quasi-polynomials of an $(\bs n, \bs \la, \bs m, \bs z)$-type.
Let $V^\star$ be the space conjugate to $V$ and defined in Section \ref{subsection conj space}.

For $a = 1, \dots , M$, let $\gamma_a$ be a small circle around $z_a$ in $\C$
oriented counterclockwise. Denote by $U$ the complex vector space
spanned by the functions of the form
\bean\label{int transform}
\hat f_a(u)\ =\ \int_{\gamma_a} x^{u}\, f(x)\, x^{-N}\prod_{a=1}^M(x-z_a)^{-1}\, dx\ ,
\eean
where $a = 1, \dots , M$, \ $f \in V^\star$.

\begin{theorem}
\label{thm special Mellin transform}
Let $(V, \bs z)$ be a space of quasi-polynomials of
an \newline
$(\bs n, (\la_1,\dots,\la_N), \bs m,
(z_1,\dots,z_M))$- type. Then

\begin{enumerate}
\item[(i)] The space $U$ is a space of quasi-exponentials of the
\newline
$(\bs m, (z_1,\dots,z_N), \bs n, (\la_1+1,\dots, \la_N+1))$-type.
\item[(ii)]
Let $\tilde D_V = \sum_{i=1}^M\sum_{j=1}^N
A_{ij}x^i (x\partial_x)^j$ be the special regularized fundamental differential
operator of $V$
where $A_{i j}$ are suitable complex numbers. Then
\bea
\sum_{i=1}^M\sum_{j=1}^N\, A_{i j}\ u^j \tau_u^i
\eea
is the special regularized fundamental difference operator of $U$.
\end{enumerate}
\end{theorem}

The proof is similar to the proof of Theorem \ref{thm Mellin transform}.

The special space of quasi-exponentials
$(U,\bs\la)$ is called {\it special bispectral dual}
to the special space of quasi-polynomials $(V,\bs z)$.

\subsection{Special Fourier-type transform}
\label{Special Fourier-type transform}
Let $(U,\bs \la)$ be a space of quasi-exponentials of an
$(\bs m, \bs z, \bs n, \bs \la)$-type.
Let $\Wrd_U$ be the discrete
Wronskian of $U$ and let $B_M(u)$ be the last coefficient of the
special regularized fundamental difference operator of $U$.

The complex vector space spanned by all functions of the form
\bean
\label{special con space}
\frac
{\tau_u( \Wrd(f_1,\dots,f_{M-1}))}
{B_M(u)\,\Wrd_U(u)}
\eean
with $f_i\in V$ has dimension $M$. This space is denoted by
$U^\bullet$ and called {\it special regularized conjugate to} $U$.

For $i = 1,\dots , N$, consider the arithmetic sequence
$$
\frak S_i\ =\ \{\,\la_i-1,\la_i,\la_i+1, \dots , \la_i+n_{i}-1\,\}
$$
consisting of $n_{i}+1$ terms.

For $i = 1, \dots , N$, fix a non-selfintersecting closed connected curve
$\gamma_i$ in $\C$ oriented counterclockwise such that it encircles
the sequence $\frak S_i$ and does not contain inside or intersect with
points of other sequences $\frak S_j$ for $j\neq i$.

Denote by $V$
the complex vector space spanned by functions of the form
$$
\hat f_i(x)\ =\ \int_{\gamma_i} x^{u}
f(u)\, du\ ,
$$
where $i = 1, \dots , n$, \ $f \in U^\bullet$. The
vector space $V$ is called {\it special bispectral dual to $U$}.

\begin{theorem}
\label{thm special F transform}
Let $(U,\bs \la)$ be a space of quasi-exponentials of
an \newline
$(\bs m, (z_1,\dots,z_M), \bs n,\allowbreak (\la_1,\dots,\la_N))$-type.
Then

\begin{enumerate}
\item[(i)] The space $U$ is a space of quasi-exponentials of the
\newline
$(\bs n, (\la_1-1,\dots, \la_N-1), \bs m,
(z_1,\dots,z_N))$-type.
\item[(ii)]
Let $\tilde D_u = \sum_{i=1}^M\sum_{j=1}^N
A_{ij}u^i \tau_u^j$ be the special regularized fundamental difference
operator of $U$
where $A_{i j}$ are suitable complex numbers. Then
\bea
\sum_{i=1}^M\sum_{j=1}^N\, A_{i j}\ x^j (x\p)^i
\eea
is the special regularized fundamental differential operator of $V$.
\end{enumerate}
\end{theorem}

The proof is similar to the proof of Theorem \ref{thm F transform}.

The special space of quasi-polynomials
$(V,\bs z)$ is called {\it special bispectral dual} to the special space of quasi-exponentials
$(U,\bs \la)$.

Theorems \ref{thm special Mellin transform} and \ref{thm special F
transform}
imply that if $V$ is a special space of quasi-polynomials and $U$ is special bispectral dual
to $V$, then $V$ is special bispectral dual to $U$. Similarly,
if $U$ is a special space of quasi-exponentials and $V$ is special bispectral dual
to $U$, then $U$ is bispectral dual to $U$.

\section{Special spaces and solutions of the Bethe ansatz equations}
\label{Special spaces and solutions of the Bethe ansatz equations}

\subsection{Critical points of master functions and special spaces of quasi-polynomials}
\label{sec master function}
Let $(V, \bs z)$ be a space of quasi-polynomials of
an $(\bs n, \bs \la, \bs m, \bs z)$-type.
We construct the associated master function
as follows.
Set
$$
\bar n_i\ =\ n_{i+1} +\dots + n_N\ ,
\qquad
i=1,\dots, N-1\ .
$$
Consider the new
$\bar n_1+\dots +\bar n_{N-1}$ auxiliary variables
$$
\bs t^{\langle \bs n\rangle} \ = \ (t^{(1)}_1, \dots , t^{(1)}_{\bar n_1},\
t^{(2)}_1, \dots , t^{(2)}_{\bar n_2},\
\dots ,\ t^{(N-1)}_1,
\dots , t^{(N-1)}_{\bar n_{N-1}})\ .
$$
Define {\it the master function}
\begin{align}
\Phi (\bs t^{\langle \bs n\rangle};\bs\la; \bs m; \bs z)\ =
\ \prod_{a=1}^M\,z_a^{m_a(\la_1+\,m_a/2)}\!
\prod_{1\le a<b\le M}\!(z_a-z_b)^{m_am_b}\;
\prod_{i=1}^N\prod_{j=1}^{\bar n_i}\ (t^{(i)}_j)^{\la_{i+1}-\la_i+1} &
\label{master function}
\\
{}\times\ \prod_{a=1}^M\,
\prod_{j=1}^{\bar n_{1}}
(t^{(1)}_j-z_a)^{-m_a}\
\prod_{i=1}^{N-1}
\prod_{j<j'}
(t^{(i)}_j-t^{(i)}_{j'})^{2}
\prod_{i=1}^{N-2}
\prod_{j=1}^{\bar n_{i}}
\prod_{j'=1}^{\bar n_{i+1}}
(t^{(i)}_j-t^{(i+1)}_{j'})^{-1} & \ .
\notag
\end{align}
The master function is symmetric with respect to the group
$\bs \Sigma_{\bs{\bar n}}=\Sigma_{\bar n_1}\times \dots\times \Sigma_{\bar
n_{N-1}}$ of permutations of variables $t^{(i)}_j$ preserving the upper index.

A point $\tn$ with complex coordinates is called {\it a critical
point} of $\Phi(\,\cdot\,;\bs \la; \bs m; \bs z)$ if the following system of
$\bar n_1+\dots+\bar n_{N-1}$ equations is satisfied
\begin{align}
\label{BAE}
\frac{\la_1-\la_2-1}{t^{(1)}_j}\ +\
\sum_{a=1}^M\frac{m_a}{t^{(1)}_j-z_a} \
-\ \sum_{j'=1,\ j'\neq j}^{\bar n_1}
\frac{2}{t_j^{(1)}-t_{j'}^{(1)}}\ +\ \sum_{j'=1}^{\bar n_2}
\frac{1}{t_j^{(1)}-t_{j'}^{(2)}}\ =\ 0\ ,
\\
\frac {\la_i-\la_{i+1}-1}{t^{(i)}_j}\
-\ \sum_{j'=1,\ j'\neq j}^{\bar n_i}\frac{2}{t_j^{(i)}-t_{j'}^{(i)}}\
+ \
\sum_{j'=1}^{\bar n_{i-1}}
\frac1{t_j^{(i)}-t_{j'}^{(i-1)}}\ +\
\sum_{j'=1}^{\bar n_{i+1}}
\frac{1}{t_j^{(i)}-t_{j'}^{(i+1)}}\ =\ 0\ ,
\notag
\\
\frac {\la_{N-1}-\la_N-1}{t^{(N-1)}_j}\ -\
\sum_{j'=1,\ j'\neq j}^{\bar n_{N-1}}
\frac{2}{t_j^{(N-1)}-t_{j'}^{(N-1)}}+\sum_{j'=1}^{\bar n_{N-2}}
\frac{1}{t_j^{(N-1)}-t_{j'}^{(N-2)}}\ =\ 0\ ,
\notag
\end{align}
where $j=1,\dots,\bar n_1$ in the first group of equations,
$i=2,\dots,N-2$ and $j=1,\dots,\bar n_i$ in the second group of
equations, $j=1,\dots, \bar n_{N-1}$ in the last group of equations.

In other words, a point $\tn$ is a critical point if
\bea
\left(
\Phi^{-1}
\frac{\partial \Phi}{\partial t_j^{(i)}}\right)(\tn; \bs \la; \bs m; \bs z)\ =\ 0\ ,
\qquad
i = 1 , \dots , N-1,\ j = 1 , \dots, \bar n_i\ .
\eea

In the Gaudin model, equations \Ref{BAE} are called {\it the Gaudin Bethe
ansatz equations}.

\medskip

The $\bs \Sigma_{\bs{ \bar n}}$-orbit of a point $\bs t^{\langle \bs n\rangle}
\in \C^{\bar n_1+\dots+\bar n_{N-1}}$ is uniquely determined by the
$N-1$-tuple $\bs y^{\bs t^{\langle \bs n\rangle}} = (y_1, \dots , y_{N-1})$ of polynomials in $x$,\
where
\bea
&
y_i\ =\ \prod_{j=1}^{\bar n_{i}}\ (x - t^{(i)}_j)\ ,
\qquad i=1,\dots, N-1\ .
\eea
We say that $\bs y$ {\it represents}
the orbit. Each polynomial of the tuple is
considered up to multiplication by a nonzero
number since we are interested in the
roots of the polynomial only.

\medskip

We say that $\bs t^{\langle \bs n\rangle} \in \C^{\bar n_1+\dots+\bar n_{N-1}}$
is {\it Gaudin admissible}
if the value $\Phi (\bs t^{\langle \bs n\rangle};\bs \la; \bs m; \bs z)$
is well defined and is not zero.

A point $\bs t^{\langle \bs n\rangle}$ is Gaudin
admissible if and only if the associated tuple has
the following properties.
\begin{enumerate}
\item[$\bullet$]
For $a=1,\dots,M$, if $m_a > 0$, then $y_1(z_a)\neq 0$.
\item[$\bullet$]
For all $i$, $y_i(0) \neq 0$.
\item[$\bullet$]
For all $i$,
the polynomial $y_i$ has no multiple roots and no common roots with $y_{i-1}$ or $y_{i+1}$.
\end{enumerate}

Such tuples are called {\it Gaudin admissible}.

\medskip

Return to $(V, \bs z)$, a space of an $(\bs n, \bs \la, \bs m, \bs z)$-type.

The space $V = \langle x^{\la_1}p_1(x), \dots , x^{\la_N}p_N(x) \rangle$
determines the $N-1$-tuple
\ $\bs y^{V} = (y_1, \dots , y_{N-1})$ of polynomials in $x$,\ {} where
\bea
y_{N-1} = p_N\,,
\qquad
y_i\ =\ x^{\frac{(N-i)(N-i-1)}2-\la_{i+1}-\dots -\la_N} \Wr\, (x^{\la_{i+1}}p_{i+1}(x),
\dots , x^{\la_{N}}p_{N}(x))
\eea
for $i=1,\dots, N-2$.

We call the special space $(V, \bs z)$ {\it Gaudin admissible} if the tuple $\bs y^V$ is
Gaudin admissible.

\begin{theorem}\cite{MV1, MV4}
\label{thm space - cr pt}
${}$

\begin{enumerate}
\item[(i)]
Assume that the special space $(V, \bs z)$ is Gaudin admissible.
Then the tuple $\bs y^V$ represents the orbit of
a critical point of the master function.
\item[(ii)]
Assume that $\bs t^{\langle \bs n \rangle}$ is Gaudin admissible
and $\bs t^{\langle \bs n\rangle}$ is a critical point of the master function. Let
$\bs y=(y_1,\dots,y_{N-1})$ be the tuple representing the orbit of $\tn$.
Then the differential operator
\begin{align*}
\bar D\ =
\ \left( \partial_x -
\ln' ( \frac { x^{\la_1-N+1}\,\prod_{a=1}^M (x-z_a)^{m_a} } { y_1 } ) \right)
\left( \partial_x - \ln' ( \frac {x^{\la_2-N+2}\,y_1 }{y_2}) \right)
\phantom{aaaaaa}
&
\\
\phantom{aaaaaa}
{}\dots\
\left( \partial_x - \ln' (\frac{ x^{\la_{N-1}-1}\,y_{N-2} }{ y_{N-1} })\right)
\left( \partial_x - \ln' (x^{\la_N} y_{N-1} ) \right) &
\end{align*}
of order $N$ is the monic fundamental differential operator of a Gaudin
admissible special space of quasi-polynomials $(V, \bs z)$ of the
$(\bs n, \bs\la, \bs m, \bs z)$-type.
\item[(iii)] The correspondence between Gaudin
admissible special spaces of quasi-polynomials of the
$(\bs n, \bs \la, \bs m, \bs z)$-type and orbits
of Gaudin admissible critical points of the master function described in parts
{\rm (i), (ii)} is reflexive.
\end{enumerate}
\end{theorem}

This theorem establishes a one-to-one correspondence between Gaudin admissible
special spaces of quasi-polynomials of the $(\bs n,\bs\la,\bs m,\bs z)$-type
and orbits of Gaudin admissible critical points of the master function.

\subsection{Solutions of the Bethe ansatz equations
and special spaces of quasi-expo- nentials}
\label{sec BAE and exponentials}
Let $(U, \bs \la)$ be a space of quasi-exponentials of
an $(\bs m, \bs z, \bs n, \bs \la)$-type.
We define the associated system of Bethe ansatz equations
as follows.
Set
$$
\bar m_a\ =\ m_{a+1} +\dots + m_M\ ,
\qquad
a=1,\dots, M-1\ .
$$
Consider the new
$\bar m_1+\dots +\bar m_{M-1}$ auxiliary variables
$$
\bs t^{\langle \bs m\rangle} \ = \ (t^{(1)}_1, \dots , t^{(1)}_{\bar m_1},\
t^{(2)}_1, \dots , t^{(2)}_{\bar m_2},\
\dots ,\ t^{(M-1)}_1,
\dots , t^{(M-1)}_{\bar m_{M-1}})\ .
$$
{\it The Bethe ansatz equations} is
the following system of
$\bar m_1+\dots+\bar m_{M-1}$ equations:
\begin{align}
\label{BAE discrete}
\prod_{i=1}^N \frac {t_b^{(1)}-\la_i}
{t_b^{(1)}-\la_i-n_i}
\prod_{b' \neq b}
\frac {t_b^{(1)} - t_{b'}^{(1)} - 1}
{t_b^{(1)} - t_{b'}^{(1)} + 1}
\prod_{b'=1}^{\bar m_{2}}
\frac{t_b^{(1)} - t_{b'}^{(2)}+1}
{t_b^{(1)} - t_{b'}^{(2)}}\ &{}=\ \frac{z_2}{z_1}\ ,
\\
\prod_{b'=1}^{\bar m_{a-1}}
\frac{t_b^{(a)} - t_{b'}^{(a-1)}}
{t_b^{(a)} - t_{b'}^{(a-1)}-1}
\prod_{b' \neq b}
\frac {t_b^{(a)} - t_{b'}^{(a)} - 1}
{t_b^{(a)} - t_{b'}^{(a)} + 1}
\prod_{b'=1}^{\bar m_{a+1}}
\frac{t_b^{(a)} - t_{b'}^{(a+1)}+1}
{t_b^{(a)} - t_{b'}^{(a+1)}}\ &{}=\ \frac{z_{a+1}}{z_a}\ ,
\notag
\\
\prod_{b'=1}^{\bar m_{N-2}}
\frac{t_b^{(N-1)} - t_{b'}^{(N-2)}}
{t_b^{(N-1)} - t_{b'}^{(N-2)}-1}
\prod_{b' \neq b}
\frac {t_b^{(N-1)} - t_{b'}^{(N-1)} - 1}
{t_b^{(N-1)} - t_{b'}^{(N-1)} + 1}\ &{}=\ \frac{z_M}{z_{M-1}}
\notag
\end{align}
where $b=1,\dots,\bar m_1$ in the first group of equations,
$a=2, \dots , M-2$ and $b=1,\dots,\bar m_a$ in the second group of
equations, $b=1,\dots, \bar m_{M-1}$ in the last group of equations.

The system of Bethe ansatz equations is symmetric with respect to the group
\linebreak
$\bs \Sigma_{\bs{ \bar m}}=\Sigma_{\bar m_1}\times \dots\times \Sigma_{\bar m_{M-1}}$
of permutations of variables $t^{(a)}_b$ preserving the upper index.

In the XXX model, equations \Ref{BAE discrete} are called {\it the Bethe
ansatz equations}.

\medskip

The $\bs \Sigma_{\bs{ \bar m}}$-orbit of a point $\bs t^{\langle \bs m\rangle}
\in \C^{\bar m_1+\dots+\bar m_{M-1}}$ is uniquely determined by the
$M-1$-tuple $\bs y^{\bs t^{\langle \bs m\rangle}}
= (y_1, \dots , y_{M-1})$ of polynomials in $u$,\
where
\bea
&
y_a\ =\ \prod_{b=1}^{\bar m_{a}}\ (u - t^{(a)}_b)\ ,
\qquad a=1,\dots, M-1\ .
\eea
We say that $\bs y$ {\it represents}
the orbit. Each polynomial of the tuple is
considered up to multiplication by a nonzero number.

\medskip

We say that $\bs t^{\langle \bs m\rangle} \in \C^{\bar m_1+\dots+\bar
m_{M-1}}$ is {\it XXX admissible} of $\bs\la$-type if
\begin{align*}
t_b^{(a)}\neq t_{b'}^{(a)}\ ,
\qquad t_b^{(a)}\neq t_{b'}^{(a)}&+1\ ,
\qquad
t_b^{(a)}\neq t_{b'}^{(a+1)}\ ,
\qquad
t_b^{(a)}\neq t_{b'}^{(a+1)}-1\ ,
\\
&t_b^{(1)}\neq \la_i + r\ ,
\end{align*}
for
all $a,b,b',i$, and $r=0,\dots,n_i$.

If $\bs t^{\langle \bs m\rangle} \in \C^{\bar m_1+\dots+\bar
m_{M-1}}$ is XXX admissible, then the corresponding tuple
$\bs y^{\bs t^{\langle \bs m\rangle}}$ is called {\it XXX admissible}.

\medskip

Return to $(U, \bs \la)$, a special space of quasi-exponentials
of an $(\bs m, \bs z, \bs n, \bs \la)$-type.

The space $U = \langle z_1^{u}q_1(u), \dots , z_M^{u}q_M(u) \rangle$
determines the $M-1$-tuple
\ $\bs y^{U} = (y_1, \dots , y_{M-1})$ of polynomials in $u$,\ {} where
\bea
y_{M-1} = q_M\, ,
\qquad
y_a\ =\ \prod_{b=a+1}^M z_b^{-u}\, \Wrd (z_{a+1}^u q_{a+1}(u), \dots, z_{M}^u q_{M}(u))
\eea
for $a = 1, \dots, M-2$.

We call the special space $(U, \bs \la)$ {\it XXX admissible}
if the tuple $\bs y^U$ is XXX admissible.

\begin{theorem}\cite{MV1, MV4}
\label{thm space - BAE}
${}$

\begin{enumerate}
\item[(i)]
Assume that the special space $(U, \bs \la)$ is XXX admissible.
Then the tuple $\bs y^U$ represents the orbit of a
solution of the XXX Bethe ansatz equations \Ref{BAE discrete}.
\item[(ii)]
Assume that $\bs t^{\langle \bs m \rangle}$ is XXX admissible and $\bs t^{\langle
\bs m\rangle}$ is a
solution of the XXX Bethe ansatz equations \ref{BAE discrete}. Let
$\bs y = (y_1,\dots,y_{M-1})$ be the tuple representing the orbit of $\bs t^{\langle
\bs m\rangle}$. Then the difference operator
\begin{align*}
\bar D\ =
\ \left( \tau_u - \frac { y_1(u) }{ y_1(u+1) }
\prod_{i=1}^N \frac {u-\la_i +1} {u-\la_i-n_i+1}\;z_1\right)
\left( \tau_u - \frac { y_1(u+1) }{ y_1(u) }
\frac { y_2(u) }{ y_2(u+1) }\;z_2\right) &
\\
{}\dots
\ \left( \tau_u - \frac { y_{M-2}(u+1) }{ y_{M-2}(u) }
\frac { y_{M-1}(u) }{ y_{M-1}(u+1) }\;z_{M-1}\right)
\left( \tau_u - \frac { y_{M-1}(u+1) }{ y_{M-1}(u)}\;z_M \right) &
\end{align*}
of order $M$ is the monic fundamental difference operator of a special
space of quasi-exponentials
$(U, \bs \la)$ of the $(\bs m, \bs z, \bs n, \bs \la)$-type.
\item[(iii)] The correspondence between
XXX admissible special spaces of quasi-exponentials of the
$( \bs m, \bs z, \bs n, \bs \la)$-type and orbits
of XXX admissible solutions of the XXX Bethe ansatz equations
described in parts
{\rm (i), (ii)} is reflexive.
\end{enumerate}
\end{theorem}

This theorem establishes a one-to-one correspondence between XXX admissible
special spaces of quasi-exponentials of the $(\bs m,\bs z,\bs n,\bs\la)$-type
and orbits of XXX admissible solutions of the XXX Bethe ansatz equations.

\section{Finiteness of solutions of the Bethe ansatz equations}
\label{Sec FINIT}

\subsection{Finiteness of the Gaudin admissible critical points}
\label{sec Finiteness of the Gaudin admissible critical points}

\begin{lem}\label{lem Gaudin finiteness}
For given fixed $\bs m,\ \bs z$ and generic $\bs \la$,
the master function $\Phi(\bs t^{\langle \bs n\rangle};\bs\la; \bs m; \bs z)$,
defined in \Ref{master function},
has only finitely many Gaudin admissible
critical points.
\end{lem}
The lemma follows from Lemma 2.1 in \cite{MV1}.

\subsection{The number of orbits of the Gaudin admissible critical points}
\label{sec number Gaidin }
Let $\la_1, \dots ,$ $\la_N\,\in\, \C$ be distinct
numbers such that $\la_i-\la_j \notin \Z$ for $i\neq j$.
Set $\bar n = n_1+\dots + n_N$. Consider the complex vector space $X$ spanned by
functions $x^{\la_i+j}$, \ $i=1, \dots , N$,
$j = 0,\dots ,n_i$. The space $X$ is of dimension
$\bar n + N$.

For $z\in\C^*$, define a complete flag $\bs F(z)$ in $X$,
\bea
\bs F(z)\ =\ \{\, 0=F_0(z) \subset F_1(z) \subset
\dots \subset F_{\bar n+N}(z)=X \,\}\ ,
\eea
where $F_k(z)$ consists of all $f\in X$ which have zero at $z$ of order not less that
$\bar n + N -k$. The subspace $F_k(z)$ has dimension $k$.

Define two complete flags of $X$ {\it at infinity}.

Say that $x^{\la_i+j}\ <_1\
x^{\la_{i'}+j'}$ if $i<i'$ or $i=i'$ and $j<j'$. Set
\bea
\bs F(\infty_1)\ =\ \{\, 0=F_0(\infty_1) \subset F_1(\infty_1) \subset
\dots \subset F_{\bar n+N}(\infty_1)=X \,\}\ ,
\eea
where $F_k(\infty_1)$ is spanned by $k$ smallest elements with respect to $<_1$.

Say that $x^{\la_i+j}\ <_2\ x^{\la_{i'}+j'}$
if $i>i'$ or $i=i'$ and $j<j'$. Set
\bea
\bs F(\infty_2)\ =\ \{ 0=F_0(\infty_2) \subset F_1(\infty_2) \subset
\dots \subset F_{\bar n+N}(\infty_2)=X \}\ ,
\eea
where $F_k(\infty_2)$ is spanned by $k$ smallest elements with respect to $<_2$.

\bigskip

Denote by $\GR$ the Grassmannian manifold of $N$-dimensional vector subspaces of $X$.
Let $\bs F$ be a complete flag of $X$,
$$
\bs F\ = \ \{\,0=F_0
\ \subset F_1\ \subset \ \dots\ \subset \ F_{\bar n +N}\ =\ X\, \}\ .
$$
{\it A ramification sequence} is a sequence
$(c_1,\dots,c_N)\in \Z^N$ such that $\bar n \geq c_1 \geq \dots \geq c_{N} \geq
0.$ For a ramification sequence $\bs c =(c_1,\dots,c_{N})$
define {\it the Schubert cell}
\begin{align}
&
\Omega^o_{\bs c}(\bs F) \ =\ \{\,V\in \GR \mid \dim(V\cap F_u) = \ell ,\
\notag
\\
&
\phantom{aaaaaaaaaaaa}
\bar n +\ell-{c}_{\ell}\ \leq\ u\ < \ \bar n + \ell + 1 - {c}_{\ell+1},
\ \text{}\ell = 0, \dots,N \}\ ,
\notag
\end{align}
where ${c}_0 = \bar n$,\ {}\ ${c}_{N+1} = 0$. The cell
$\Omega^o_{\bs c}(\bs F)$ is a smooth connected variety.
The closure of $\Omega^o_{\bs c}(\bs F)$ is denoted by
$\Omega_{\bs c}(\bs F)$. The codimension
of $\Omega^o_{\bs c}(\bs F)$ is
$$
|\bs c|\ =\ {c}_1\ +\
{c}_2\ +\ \dots\ +\ {c}_{N}\ .
$$
Every
$N$-dimensional vector subspace of $X$ belongs to a unique Schubert cell
$\Omega^o_{\bs c}(\bs F)$.

\bigskip

For $a=1,\dots , M$, define the ramification sequence
\bea
\bs c(a)\ =\ ( m_a , 0 , \dots , 0 )\ .
\eea
Define the ramification sequences
\bea
\bs c(\infty_1)& =& (n_2+\dots+n_N, n_3+\dots+ n_N,\dots,n_{N},0)\ ,
\\
\bs c(\infty_2) & = & (n_1+\dots+n_{N-1}, n_1+\dots+ n_{N-2},\dots,n_{1},0)\ .
\eea

\begin{lem} ${}$

\begin{enumerate}
\item[$\bullet$] We have
\bea
&&
\sum_{a=1}^M \ {\rm codim}\ \Omega^o_{\bs c(a)}(\bs F(z_a))\ +\
{\rm codim}\ \Omega^o_{\bs c(\infty_1)}(\bs F(\infty_1))
+
\\
&&
\phantom{aaaaaaaaaaaaaa}
{\rm codim}\ \Omega^o_{\bs c(\infty_2)}(\bs F(\infty_2))
=
\dim\ \GR\ =\ N \bar n\ .
\eea
\item[$\bullet$]
Let $V\in\GR$. The pair $(V, \bs z)$ is a space of
the $(\bs n, \bs \la, \bs m, \bs z)$-type, if and only if $V$ belongs
to the intersection of $M+2$ Schubert cells \bea && \Omega^o_{\bs
c(1)}(\bs F(z_1)) \cap \Omega^o_{\bs c(2)}(\bs F(z_2)) \cap \dots
\phantom{aaaaaaaaaaaaaa} \\ && \phantom{aaaaaaaa} \cap \Omega^o_{\bs
c(M)}(\bs F(z_M)) \cap \Omega^o_{\bs c(\infty_1)}(\bs F(\infty_1))
\cap \Omega^o_{\bs c(\infty_2)}(\bs F(\infty_2)) \ . \eea
\end{enumerate}
\hfill
$\square$
\end{lem}

According to Schubert calculus,
the multiplicity of the intersection of Schubert cycles
\bean\label{intersection}
&&
\Omega_{\bs c(1)}(\bs F(z_1)) \cap \Omega_{\bs c(2)}(\bs F(z_2)) \cap
\dots
\notag
\\
&&
\phantom{aaaaaa}
\cap
\Omega_{\bs c(M)}(\bs F(z_M)) \cap
\Omega_{\bs c(\infty_1)}(\bs F(\infty_1)) \cap
\Omega_{\bs c(\infty_2)}(\bs F(\infty_2)) \
\eean
can be expressed in representation-theoretic terms as follows.

For a ramification sequence $\bs c$ denote by $L^{\langle N \rangle}_{\bs c}$ the finite
dimensional irreducible \linebreak
$\glN$-module with highest weight $\bs c$.
Any $\glN$-module $L^{\langle N \rangle}$ has a natural structure of an $\slN$-module
denoted by $\widetilde L^{\langle N \rangle}$. By \cite{Fu},
the multiplicity of the intersection in
\Ref{intersection} is equal to the multiplicity of the trivial
$\slN$-module
in the tensor product of $\slN$-modules
\bean\label{slN intersection}
\widetilde L^{\langle N \rangle}_{\bs c(1)} \ox\dots \ox\widetilde L^{\langle N \rangle}
_{\bs c(M)} \ox
\widetilde L^{\langle N \rangle}_{\bs c(\infty_1)} \ox \widetilde
L^{\langle N \rangle}_{\bs c(\infty_2)} \ .
\eean

\begin{prop} [\cite{MTV2}]
\label{lem on multiplicity}
The multiplicity of the trivial $\slN$-module in the tensor product
\Ref{slN intersection} is equal to the dimension of the weight
subspace of weight $[n_1, \dots , n_N]$ in the tensor product of
$\glN$-modules
\bean
\label{glN intersection}
L^{\langle N \rangle}_{\bs c(1)} \ox\dots \ox L^{\langle N \rangle}_{\bs c(M)}\ .
\eean
\end{prop}

\begin{cor}\label{cor bound}
For generic $\bs \la$, the number of orbits of
the Gaudin admissible critical points of
the master function
$\Phi(\bs t^{\langle \bs n\rangle};\bs\la; \bs m; \bs z)$,
is not greater than
the dimension of the weight space
$ (L^{\langle N \rangle}_{\bs c(1)}
\ox\dots \ox L^{\langle N \rangle}_{\bs c(M)})[n_1,\dots,n_N]$.
\hfill $\square$
\end{cor}

\subsection{Finiteness of solutions of Bethe ansatz equations
\Ref{BAE discrete}}
\label{sec Finiteness of the XXX }

Let $\bs z = (z_1,\dots, z_M)$ $ \in \C^M$. Let $(\bar m_1,\dots, \bar m_{M-1})$
be a collection of nonnegative integers. We say that
$\bs z$ {\it is separating with respect to } $(\bar m_1,\dots, \bar m_{M-1})$
if
$$
\prod_{a=1}^{M-1}\ \left( \frac{z_{a+1}}{z_a}\right)^{c_a}\ \neq \ 1
$$
for all sets of integers $\{c_1, \dots , c_{M-1} \}$ such that
$0 \leq c_a \leq \bar m_a$, $\sum_a c_a > 0$.

Clearly, for given $(\bar m_1,\dots, \bar m_{M-1})$, a generic $\bs z$
is separating.

\begin{lem}
Let $\bs n$, $\bs\la$,\ $(\bar m_1,\dots, \bar m_{M-1})$ be fixed.
Let $\bs z$ be separating with respect to
$(\bar m_1,\dots, \bar m_{M-1})$. Then
the Bethe ansatz equations \Ref{BAE discrete}
have only finitely many XXX admissible solutions of $\bs \la$-type.
\end{lem}

\begin{proof}
If the algebraic set of XXX admissible solutions of
\Ref{BAE discrete}
is infinite, then it is unbounded.
Suppose that we have a sequence of solutions which is
unbounded. Without loss of generality, we assume that $t_{b}^{(a)}$ tends
to infinity for $a = 1, \dots , M-1, \ b = 1, \dots , c_a$,
and remains bounded for
all other values of $a, b$.
Multiply all of the equations in \Ref{BAE discrete}
corresponding to a variables $t_b^{(a)}$ with
$a = 1, \dots , M-1, \ b = 1, \dots , c_a$, and pass in the product
to the limit along our sequence of solutions.
Then the resulting equation is
$$
\prod_{a=1}^{M-1}\ \left( \frac{z_{a+1}}{z_a}\right)^{c_a}\ = \ 1\ .
$$
This equation contradicts to our assumption.
\end{proof}

\subsection{The number of orbits of the XXX admissible solutions}
\label{sec number XXX }
Let $z_1, \dots ,$ $z_M\,\in\, \C$ be distinct
numbers with fixed argument.
Set $\bar m = m_1+\dots + m_M$. Consider the complex vector space $Y$ spanned by
functions $z_a^uu^{b}$, \ $a=1, \dots , M$,
$b = 0,\dots ,m_a$. The space $Y$ is of dimension
$\bar m + M$.

For $\la\in\C$, define a complete flag $\bs F(\la)$ in $Y$,
\bea
\bs F(\la)\ =\ \{\, 0=F_0(\la) \subset F_1(\la) \subset
\dots \subset F_{\bar m+M}(\la)=Y \,\}\ ,
\eea
where $F_k(\la)$ consists of all $f\in Y$ which are divisible by
\ $\prod_{j=1}^{\bar m + M - k}\,(u-\la - j+ 1)$.
The subspace $F_k(\la)$ has dimension $k$.

Define two complete flags of $Y$ {\it at infinity}.

Say that $z_a^uu^{b}\ <_1\ z_{a'}^uu^{b'}$
if $a<a'$ or $a=a'$ and $b<b'$. Set
\bea
\bs F(\infty_1)\ =\ \{\, 0 = F_0(\infty_1)
\subset F_1(\infty_1) \subset
\dots \subset F_{\bar m+M}(\infty_1)=Y \,\}\ ,
\eea
where $F_k(\infty_1)$ is spanned by $k$ smallest elements with respect to $<_1$.

Say that $z_a^uu^{b}\ <_2\ z_{a'}^uu^{b'}$
if $a>a'$ or $a=a'$ and $b<b'$. Set
\bea
\bs F(\infty_2)\ =\ \{\, 0=F_0(\infty_2) \subset F_1(\infty_2) \subset
\dots \subset F_{\bar m + M}(\infty_2)=Y \,\}\ ,
\eea
where $F_k(\infty_2)$ is spanned by $k$ smallest elements with respect to $<_2$.

\bigskip

Denote by $\GRM$ the Grassmannian manifold of $M$-dimensional vector subspaces of $Y$.
Let $\bs F$ be a complete flag of $Y$,
$$
\bs F\ = \ \{\,0=F_0
\ \subset F_1\ \subset \ \dots\ \subset \ F_{\bar m +M}\ =\ Y\, \}\ .
$$
For a ramification sequence $\bs c =(c_1,\dots,c_{M}) \in \Z^M$,
$\bar m \geq c_1 \geq \dots \geq c_{M} \geq 0,$ denote by $\Omega^o_{\bs c}(\bs F)
\subset \GRM$ the corresponding Schubert cell.

\bigskip

For $i=1,\dots , N$, define the ramification sequence
\bea
\bs c(i)\ =\ ( n_i , 0 , \dots , 0 )\ .
\eea
Define the ramification sequences
\bea
\bs c(\infty_1)& =& (m_2+\dots+m_M, m_3+\dots+ m_M,\dots,m_{M},0)\ ,
\\
\bs c(\infty_2) & = & (m_1+\dots+m_{M-1}, m_1+\dots+ m_{M-2},\dots, m_{1},0)\ .
\eea

\begin{lem} ${}$

\begin{enumerate}
\item[$\bullet$] We have
\bea
&&
\sum_{i=1}^N \ {\rm codim}\ \Omega^o_{\bs c(i)}(\bs F(\la_i))\ +\
{\rm codim}\ \Omega^o_{\bs c(\infty_1)}(\bs F(\infty_1))
+
\\
&&
\phantom{aaaaaaaaaaaaaa}
{\rm codim}\ \Omega^o_{\bs c(\infty_2)}(\bs F(\infty_2))
=
\dim\ \GRM\ =\ M \bar m\ .
\eea
\item[$\bullet$]
Let $U\in\GRM$. The pair $(U, \bs \la)$ is a space of
the $(\bs m, \bs z, \bs n, \bs \la)$-type, if and only if $U$ belongs
to the intersection of $N+2$ Schubert cells
\bea
&& \Omega^o_{\bs
c(1)}(\bs F(\la_1)) \cap \Omega^o_{\bs c(2)}(\bs F(\la_2)) \cap \dots
\phantom{aaaaaaaaaaaaaa}
\\
&&
\phantom{aaaaaaaa}
\cap \Omega^o_{\bs
c(M)}(\bs F(\la_N)) \cap \Omega^o_{\bs c(\infty_1)}(\bs F(\infty_1))
\cap \Omega^o_{\bs c(\infty_2)}(\bs F(\infty_2)) \ .
\eea
\end{enumerate}
\hfill
$\square$
\end{lem}

According to Schubert calculus,
the multiplicity of the intersection of Schubert cycles
\bean
\label{XXX intersection}
&& \Omega_{\bs
c(1)}(\bs F(\la_1)) \cap \Omega_{\bs c(2)}(\bs F(\la_2)) \cap \dots
\phantom{aaaaaaaaaaaaaa}
\notag
\\
&&
\phantom{aaaaaaaa}
\cap \Omega_{\bs
c(M)}(\bs F(\la_N)) \cap \Omega_{\bs c(\infty_1)}(\bs F(\infty_1))
\cap \Omega_{\bs c(\infty_2)}(\bs F(\infty_2)) \ .
\eean
can be expressed in representation-theoretic terms as follows.

For a ramification sequence $\bs c = (c_1,\dots,c_{M}) \in \Z^M$,
denote by $L^{\langle M \rangle}_{\bs c}$ the finite
dimensional irreducible
$\glM$-module with highest weight $\bs c$.
Any $\glM$-module $L^{\langle M \rangle}$ has a natural structure of an $\slM$-module
denoted by $\widetilde L^{\langle M \rangle}$. By \cite{Fu},
the multiplicity of the intersection in
\Ref{intersection} is equal to the multiplicity of the trivial
$\slM$-module
in the tensor product of $\slM$-modules
\bean\label{slM intersection}
\widetilde L^{\langle M \rangle}_{\bs c(1)}
\ox\dots \ox\widetilde L^{\langle M \rangle}_{\bs c(N)} \ox
\widetilde L^{\langle M \rangle}_{\bs c(\infty_1)}
\ox \widetilde L^{\langle M \rangle}_{\bs c(\infty_2)} \ .
\eean
By Proposition \ref{lem on multiplicity},
the multiplicity of the trivial $\slM$-module in the tensor product
\Ref{slM intersection} is equal to the dimension of the weight
subspace of weight $[m_1, \dots , m_M]$ in the tensor product of
$\glM$-modules
\bean
\label{glM intersection}
L^{\langle M \rangle}_{\bs c(1)} \ox\dots \ox L^{\langle M \rangle}_{\bs c(N)}\ .
\eean

\begin{cor}\label{cor XXX bound}
For generic $\bs z$, the number of orbits of
the XXX admissible of $\bs \la$-type solutions of
the Bethe ansatz equations \Ref {BAE discrete}
is not greater than
the dimension of the weight space
$ (L^{\langle M \rangle}_{\bs c(1)}
\ox\dots \ox L^{\langle M \rangle}_{\bs c(N)})[m_1,\dots,m_M]$.
\hfill $\square$
\end{cor}

\section{The KZ and dynamical Hamiltonians}
\label{section KZ hamiltonians}

\subsection{The Gaudin KZ Hamiltonians}
\label{sec Gaudin KZ Hamiltonians}

Let $E_{ij}$, $i,j=1,\dots , N$, be the standard generators of the
complex Lie algebra $\frak {gl}_N$.

We have the root decomposition \,$\glN=\n^+\oplus\h\oplus\n^-$ where
$$
\n^+\>=\,\oplus _{i<j}\ \C\cdot E_{i j}\,,\qquad
\h\,=\,\oplus_{i=1}^N \ \C \cdot E_{i i}\,,\qquad
\n^-\>=\,\oplus _{i>j}\ \C \cdot E_{i j}\,.
$$
Let
\bea
\Omega^0 =
\frac 12 \sum_{i=1}^N E_{ii}\ox E_{ii}\ ,
\qquad
\Omega^+ = \Omega ^0 + \sum_{i<j} E_{ij} \ox E_{ji}\ ,
\qquad
\Omega^- = \Omega ^0 + \sum_{i<j} E_{ji} \ox E_{ij}\ .
\eea
Let $Y = Y_1\ox\dots\ox Y_M$ be the tensor product of
finite-dimensional irreducible $\glN$-modules.

{\it The Gaudin KZ Hamiltonians}
$H^{\frak G}_a(\bs\la,\bs z)$, $a=1,\dots, M$, acting on $Y$-valued
functions of \ $\bs \la=(\la_1,\dots , \la_N), \ \bs z=(z_1,\dots, z_M) $ \ are defined
by the formula \ \cite{TV5}\ :
$$
H^{\frak G}_a(\bs\la, \bs z)\,=\
\sum_{i=1}^N (\la_i - \frac {E_{ii}}2) (E_{ii})^{(a)} +
\sum_{b=1,\ b \ne a}^M
\frac{z_a (\Om^+)^{(ab)} + z_b (\Om^-)^{(a b)}}
{z_a - z_b} \ .
$$
Here the linear operator
$(\Om^\pm)^{(a b)} : Y \to Y$ acts as $\Om^\pm$ on $Y_a\ox Y_b$,
and as the identity on other tensor factors of $Y$. Similarly,
$E_{i i}^{(a)}$ acts as $E_{i i}$ on $Y_a$ and as the identity on other
factors.

\subsection{The Gaudin dynamical Hamiltonians}
\label{sec Gaudin dynamical Hamiltonians}
For any $i,j = 1,\dots,N,\ i\neq j$, introduce a series $B_{ij}(t)$
depending on a complex number $t$:
\bea
B_{i,j}(t)\ =\ 1\ +\ \sum_{s=1}^\infty\ (E_{ji})^s(E_{ij})^s\
\prod_{l=1}^s\ \frac 1
{j\,(t-E_{ii} +E_{jj} - l)}\ .
\eea
The series has a well-defined action in any finite-dimensional $\glN$-module $W$ giving
an ${\rm End}\,(W)$-valued rational function of $t$.

{\it The Gaudin
dynamical Hamiltonians} $G^{\frak G}_i(\bs\la, \bs z)$,
$i=1,\dots,N$, acting on $Y$-valued functions of $\bs \la, \bs z$ are defined by the formula
\ \cite{TV5}\ :
\bea
G_i^{\frak G}(\bs \la,\bs z) &=&
(B_{i,N}(\la_i-\la_N)\dots B_{i,\,i+1}(\la_i-\la_{i+1}))^{-1} \times
\\
&&
\phantom{aaaa}
\prod_{a=1}^M\ (z_a^{-E_{ii}})^{(a)}\ \times \
B_{1,i}(\la_1-\la_i)\dots B_{i-1,\,i}(\la_{i-1}-\la_{i}) \ .
\eea

\subsection{The Gaudin diagonalization problem}
\label{sec Gaudin diagonalization problem}
The Gaudin KZ and dynamical Hamiltonians commute \cite{TV4},
$$
[H^{\frak G}_a(\bs\la,\bs z) , H_b^{\frak G}(\bs \la, \bs z)]\ =\ 0\ ,
\quad [H_a^{\frak G}(\bs\la,\bs z)\>,G_i^{\frak G}(\bs\la,\bs z)]\ =\ 0\ ,
\quad
[G_i^{\frak G}(\bs \la,\bs z) ,G_j^{\frak G}(\bs\la,\bs z)]\ =\ 0\ ,
$$
for $a, b = 1,\dots, M$, and $i, j = 1,\dots,N$.

{\it The Gaudin diagonalization problem} is to diagonalize simultaneously
the Gaudin KZ Hamiltonians $H_a^{\frak G}$, $a=1,\dots, M$,\
and dynamical Hamiltonians $ G_i^{\frak G}$,
$i=1,\dots, N$, for given $\bs\la, \bs z$.
The Hamiltonians preserve the weight decomposition of $Y$
and the diagonalization problem can be considered on a given weight subspace of $Y$.

\subsection{Diagonalization and critical points}
For a nonnegative integer $m$, denote by $L_m^{\langle N\rangle}$ the irreducible
finite-dimensional
$\glN$-module with highest weight $(m,0,\dots,0)$.

Let $\bs n=(n_1,\dots,n_N)$ and $\bs m=(m_1,\dots,m_M)$ be vectors of nonnegative integers
with $\sum_{i=1}^N n_i = \sum_{a=1}^M m_a$.

Consider the tensor product $L_{m_1}^{\langle N\rangle}\ox\dots\ox
L_{m_M}^{\langle N\rangle}$ and its weight subspace
\bea
\LMN\ =\ (L_{m_1}^{\langle N\rangle}\ox\dots\ox L_{m_M}^{\langle N\rangle})
[n_1,\dots,n_N]\ .
\eea
Let $\bs \la \in \C^N$, $\bs z\in\C^M$. Assume that each of $\bs\la$ and
$ \bs z$ has distinct coordinates.
Consider the Hamiltonians $H_a^{\frak G}
(\bs \la,\bs z)$,
$G_i^{\frak G}(\bs \la,\bs z)$
acting on $\LMN$.
{\it The Bethe ansatz method} is a method
to construct common eigenvectors of the Hamiltonians.

As in Section \ref{sec master function} consider the space $\C^{\bar
n_1+\dots +\bar n_{N-1}}$ with coordinates $\bs t^{\langle \bs
n\rangle}$. Let
\linebreak
$\Phi(\bs t^{\langle \bs
n\rangle};\bs\la;\bs m;\bs z)$ be the master
function on $\C^{\bar n_1+\dots +\bar n_{N-1}}$ defined in \Ref{master
function}.

In Section 4 of \cite{MaV}, a certain $\LMN$-valued
rational function
$\omega^{\frak G} : \C^{\bar n_1+\dots +\bar n_{N-1}}\to \LMN$,
depending on $\bs z$,
is constructed. It is called {\it the Gaudin universal rational function}.
For $N=2$ formulas for the Gaudin universal rational
function see below in Section \ref{gaudin rational function}.

The Gaudin universal rational function is well defined for admissible
$\bs t^{\langle \bs n\rangle}$. The Gaudin universal
rational function is symmetric with respect to
the $\bs \Sigma_{\bs{ \bar n}}$-action.

\begin{theorem}[ \cite{RV, MaV} ]
\label{thm RV}
If $\bs t^{\langle \bs n\rangle}$ is a Gaudin admissible
non-degenerate critical
point of the master function $\Phi(\cdot\,;\bs\la-\bs n;\bs m;\bs z)$,
then $\omega^{\frak G} (\bs t^{\langle \bs n\rangle},\bs z)\in \LMN$ is
an eigenvector of the Gaudin Hamiltonians,
\begin{alignat}2
& H_a^{\frak G}(\bs \la,\bs z)\ \omega^{\frak G}(\bs t^{\langle \bs n\rangle},
\bs z)\ =\ z_a\,\frac{\partial }{\partial z_a}\ {\rm log}\
\Phi (\bs t^{\langle \bs n\rangle}; \bs\la-\bs n;\bs m;\bs z)\
\omega^{\frak G}( && \bs t^{\langle \bs n\rangle}, \bs z)\ ,
\\
&&& a=1,\dots , M\ ,
\notag\\
& G_1^{\frak G}(\bs \la,\bs z) \ \omega^{\frak G}(\bs t^{\langle \bs n\rangle},
\bs z)\ =\ \prod_{a=1}^M\,z_a^{-m_a}\,\prod_{j=1}^{\bar n_1}\,t^{(i)}_j
\ \omega^{\frak G}(\bs t^{\langle \bs n\rangle}, \bs z)\ ,
\notag
\\
& G_i^{\frak G}(\bs \la,\bs z) \ \omega^{\frak G}(\bs t^{\langle \bs n\rangle},
\bs z)\ =\ \prod_{j=1}^{\bar n_{i-1}}\,(t^{(i-1)}_j)^{-1}\,
\prod_{j=1}^{\bar n_i}\,t^{(1)}_j
\ \omega^{\frak G}(\bs t^{\langle \bs n\rangle}, \bs z)\ ,&& i=2,\dots , N-1\ .
\notag
\\
& G_N^{\frak G}(\bs \la,\bs z) \ \omega^{\frak G}(\bs t^{\langle \bs n\rangle},
\bs z)\ =\ \prod_{j=1}^{\bar n_{N-1}}\,(t^{(N-1)}_j)^{-1}
\ \omega^{\frak G}(\bs t^{\langle \bs n\rangle}, \bs z)\ .
\notag
\end{alignat}
\end{theorem}

If $\bs t^{\langle \bs n\rangle}$ is a Gaudin admissible
critical point of the master function, then the vector
$\omega^{\frak G}(\bs t^{\langle \bs n\rangle}, \bs z)$ is called
{\it a Gaudin Bethe eigenvector}.

By Corollary \ref{cor bound}, for generic $\bs \la$, the number of Gaudin Bethe
eigenvectors is not greater than the dimension of the space $\LMN$. The Bethe
ansatz conjecture says that all eigenvectors of the Gaudin Hamiltonians are the
Gaudin Bethe eigenvectors in an appropriate sense, cf.~\cite{MV2}.

\subsection{The XXX KZ Hamiltonians}
\label{sec XXX KZ Hamiltonians}

Let $E_{ab}$, $a,b=1,\dots , M$, be the standard generators of the
complex Lie algebra $\frak {gl}_M$.

We have the root decomposition \,$\glM=\n^+\oplus\h\oplus\n^-$ where
$$
\n^+\>=\,\oplus _{a<b}\ \C\cdot E_{ab}\,,\qquad
\h\,=\,\oplus_{a=1}^M \ \C \cdot E_{aa}\,,\qquad
\n^-\>=\,\oplus _{a>b}\ \C \cdot E_{ab}\,.
$$

Let $V, W$ be irreducible finite-dimensional
$\glM$-modules with highest weight vectors
$v \in V,\ w\in W$. The associated rational R-matrix is a rational
${\rm End}\, (V\ox W)$-valued function $R_{V W}(t)$
of a complex variable $t$ uniquely determined by the
$\glM$-invariance condition,
\bea
[\,R_{V W}(t)\,,\, g\ox 1 + 1\ox g\,]\ =\ 0
\qquad
{\rm for\ any}
\quad
g \in \glM\ ,
\eea
the commutation relations
\bea
R_{V W}(t)
\ (\,t\, E_{ab}\ox 1 + \sum_{c=1}^M E_{ac}\ox E_{cb}\,)\
=\
(\,t\, E_{ij}\ox 1 + \sum_{c=1}^M E_{cb}\ox E_{ac}\,)\
R_{V W}(t)\ ,
\eea
for any $a, b$, \ and the normalization condition
\bea
R_{V W}(t) \ v\ox w \ =\ v\ox w\ .
\eea

Let $Y = Y_1\ox\dots\ox Y_N$ be the tensor product of
finite-dimensional irreducible $\glM$-modules.

{\it The XXX KZ Hamiltonians}
$H_i^{\frak X}(\bs z,\bs \la)$, $i=1,\dots, N$, acting on $Y$-valued
functions of \ $\bs z=(z_1,\dots, z_M), \ \bs \la=(\la_1,\dots , \la_N)$
\ are defined by the formula \ \cite{TV5}\ :
\begin{align*}
H_i^{\frak X}(\bs z, \bs \la)\ =\ {}&
(R_{i,N}(\la_i-\la_N)\dots R_{i,\,i+1}(\la_i-\la_{i+1}))^{-1}\,\times{}
\\
& \prod_{a=1}^M\ (z_a^{-E_{aa}})^{(i)}\times
R_{1,i}(\la_1-\la_i)\dots R_{i-1,\,i}(\la_{i-1}-\la_{i}) \ ,
\end{align*}
where $R_{jk}(t)$ is the endomorphisms of $Y$ which acts as $R_{Y_j Y_k}(t)$
on $Y_j \ox Y_k$ and as the identity on other factors of $Y$.

\subsection{The XXX dynamical Hamiltonians}
\label{sec XXX dynamical Hamiltonians}

{\it The XXX dynamical Hamiltonians}
\linebreak
$G_a^{\frak X}(\bs z, \bs \la)$, $a=1,\dots,M$, acting on $Y$-valued functions
of $\bs z, \bs \la$ are defined by the formula \ \cite{TV5}\ :
\bea
G_a^{\frak X}(\bs z,\bs \la) &=& - \frac {E_{aa}^2}2\;+\;
\sum_{i=1}^N\ \la_i\, E_{aa}^{(i)}\; +{}
\\
&&
\sum_{b=1,\ b \ne a}^M
\frac {z_b}{z_a -z_b} (E_{ab} E_{ba}- E_{aa})\;+\;
\sum_{b=1}^M \sum_{i<j}\ (E_{ab})^{(i)}(E_{ba})^{(j)}\ .
\eea

\subsection{The XXX diagonalization problem}
\label{sec XXX diagonalization problem}
The XXX KZ and dynamical Hamiltonians commute \cite{TV6},
$$
[H^{\frak X}_i(\bs z,\bs \la) , H^{\frak X}_j(\bs z,\bs \la)]\ =\ 0\ ,
\quad
[H_i^{\frak X}(\bs z,\bs \la)\>,G_a^{\frak X}(\bs z,\bs \la)]\ =\ 0\ ,
\quad
[G_a^{\frak X}(\bs z,\bs \la) ,G_b^{\frak X}(\bs z,\bs \la)]\ =\ 0\ ,
$$
for $i, j = 1,\dots,N$, and $a, b = 1,\dots, M$.

{\it The XXX diagonalization problem} is to diagonalize simultaneously
the XXX Hamiltonians for given $\bs z, \bs \la$.
The Hamiltonians preserve the weight decomposition of $Y$
and the diagonalization problem can be considered on a given weight subspace of $Y$.

\subsection{Diagonalization and solutions of the Bethe ansatz equations}
For a nonnegative integer $n$, denote by $L_n^{\langle M\rangle}$ the irreducible
$\glM$-module with highest weight $(n,0,\dots,0)$.

Let $\bs m=(m_1,\dots,m_M)$ and $\bs n=(n_1,\dots,n_N)$ be vectors of nonnegative integers
with $\sum_{a=1}^M m_a = \sum_{i=1}^N n_i$.

Consider the tensor product $L_{n_1}^{\langle M\rangle}\ox\dots\ox
L_{n_N}^{\langle M\rangle}$ and its weight subspace
\bea
\LNM\ =\ (L_{n_1}^{\langle M\rangle}\ox\dots\ox L_{n_N}^{\langle M\rangle})
[m_1,\dots,m_M]\ .
\eea
Let $\bs z\in\C^M$, $\bs \la \in \C^N$. Assume that each of $\bs z$ and
$ \bs \la$ has distinct coordinates.
Consider the Hamiltonians $H_i^{\frak X}(\bs z,\bs \la)$,
$G_a^{\frak X}(\bs z,\bs \la)$
acting on $\LNM$.
{\it The Bethe ansatz method} is a method
to construct common eigenvectors of the Hamiltonians.

As in Section \ref{sec BAE and exponentials} consider the space
$\C^{\bar m_1+\dots +\bar m_{M-1}}$ with coordinates $\bs t^{\langle \bs m\rangle}$
and the XXX Bethe ansatz equations \Ref{BAE discrete}.
Let $\;\bs1\,=\,(1,\dots,1)\,$ and
\beq
\label{xi}
\xi^{\langle \bs m\rangle} \ = \ (\xi^{(1)}_1, \dots ,
\xi^{(1)}_{\bar m_1},\ \xi^{(2)}_1, \dots , \xi^{(2)}_{\bar m_2},\
\dots ,\ \xi^{(M-1)}_1,
\dots , \xi^{(M-1)}_{\bar m_{M-1}})\ .
\eeq
be the point with coordinates $\xi^{(i)}_a=i$, for $i=1,\dots, M-1$
and $a=1,\dots, \bar m_{M-1}$.

In \cite{TV1}, a certain $\LNM$-valued rational function
$\omega^{\frak X} : \C^{\bar m_1+\dots +\bar m_{M-1}}\to \LNM\ $,
depending on $\bs\la$\,, is constructed. It is called {\it the XXX universal
rational function}. For $N=2$ formulas for the XXX universal rational function
see below in Section \ref{XXX rational function}.

The XXX universal rational function
$\omega^{\frak X}(\bs t^{\langle \bs m \rangle},\bs\la)$ is well defined if
$\bs t^{\langle \bs m \rangle}\!+\xi^{\langle \bs m \rangle}$ is XXX admissible
of $(\bs\la-\bs n+\bs1)$-type.
The XXX universal rational function is symmetric with respect to
the $\bs \Sigma_{\bs{ \bar m}}$-action.

\begin{theorem} \cite{TV2, MTV3}
\label{thm XXX bethe}
If $\bs t^{\langle \bs m\rangle}\!+\xi^{\langle \bs m \rangle}$ is an XXX
admissible solution of $(\bs\la-\bs n+\bs1)$-type of the Bethe ansatz equations
\Ref{BAE discrete}, then
$\omega^{\frak X}(\bs t^{\langle\bs m \rangle},\bs\la)\in\LNM$
is an eigenvector of the XXX Hamiltonians,
\begin{alignat}2
\label{XXX eigenvalue}
& H_i^{\frak X}(\bs z,\bs\la)
\ \omega^{\frak X}(\bs t^{\langle \bs m \rangle},\bs\la)\ =
\ \prod_{j=1}^{\bar m_1}\,\frac{t^{(1)}_j-\la_i}{t^{(1)}_j-\la_i+n_i}
\ \omega^{\frak X}(\bs t^{\langle \bs m\rangle},\bs\la)\ ,\qquad &&
i=1,\dots, N\ ,
\\[8pt]
& G_a^{\frak X}(\bs z,\bs \la)
\ \omega^{\frak X}(\bs t^{\langle \bs m \rangle},\bs\la)\ =
d_a(\bs t^{\langle \bs m\rangle},\bs\la)\
\omega^{\frak X}(\bs t^{\langle \bs m\rangle},\bs\la)\ , && a = 1,\dots , M\ ,
\notag
\end{alignat}
where
\begin{gather}
d_1(\bs t^{\langle \bs m\rangle},\bs\la)\ =
\ \sum_{i=1}^N\,n_i\Bigl(\la_i-\frac{n_i}2\Bigr)\;+\;
\sum_{j=1}^{\bar m_1}\,t^{(1)}_j\;-\;\frac{\bar m_1}2\;-
\,\sum_{b=2}^M\,\frac{m_bz_b}{z_1-z_b}\ ,
\notag
\\
\begin{aligned}
d_a(\bs t^{\langle \bs m\rangle},\bs\la)\ =
\ \sum_{j=1}^{\bar m_a}\,t^{(a)}_j\;-\,\sum_{j=1}^{\bar m_{a-1}}\,t^{(a-1)}_j
\;-\;\frac{\bar m_{a-1}+\bar m_a}2\;-\,
\sum_{b=1}^{a-1}\,\frac{m_az_b}{z_a-z_b}\;-\,
\sum_{b=a+1}^M\,\frac{m_bz_b}{z_a-z_b}\ , &
\\[6pt]
a=2,\dots, M-1\,, &
\end{aligned}
\notag
\\
d_M(\bs t^{\langle \bs m\rangle},\bs\la)\ =
\;{}-\,\sum_{j=1}^{\bar m_{M-1}}\,t^{(a-1)}_j\;-\;\frac{\bar m_{M-1}}2\;-\,
\sum_{b=1}^{M-1}\,\frac{m_M z_b}{z_a-z_b}\ .
\notag
\end{gather}
\end{theorem}

If $\bs t^{\langle \bs m\rangle}\!+\xi^{\langle \bs m \rangle}$
is an XXX admissible solution of $(\bs\la-\bs n+\bs1)$-type
of the XXX Bethe ansatz equations, then the vector
$\omega^{\frak X}(\bs t^{\langle \bs m\rangle},\bs\la)$
is called {\it an XXX Bethe eigenvector}.

By Corollary \ref{cor XXX bound}, for generic $\bs z$
the number of the XXX Bethe eigenvectors is
not greater than the dimension of the space $\LNM$. The Bethe ansatz conjecture
says that all eigenvectors of the XXX Hamiltonians are the XXX Bethe
eigenvectors.

\bigskip

According to previous discussions, eigenvectors of Gaudin Hamiltonians
on $\LMN$ are related to critical points of the master function
$\Phi(\bs t^{\langle \bs n\rangle};\bs\la-\bs n;\bs m;\bs z)$. The critical points
of $\Phi(\bs t^{\langle \bs n\rangle};\bs\la-\bs n;\bs m;\bs z)$ are related
to spaces of the $(\bs n, \bs\la-\bs n, \bs m, \bs z)$-type. The spaces
of the $(\bs n,\bs\la-\bs n,\bs m,\bs z)$-type are special bispectral dual
to spaces of the $(\bs m, \bs z, \bs n, \bs\la-\bs n+\bs 1)$-type, where
$\;\bs1=(1,\dots, 1)\in\C^N\;$. The spaces of the
$(\bs m,\bs z,\bs n,\bs\la-\bs n+\bs 1)$-type are related to solutions
of the XXX Bethe ansatz equations, which in its turn are related to
eigenvectors of the XXX Hamiltonians acting on $\LNM$. As a result of this
chain of relations, the eigenvectors of the Gaudin Hamiltonians on $\LMN$ and
the eigenvectors of the XXX Hamiltonians on $\LNM$ must be related.

Indeed this relation is given by the $(\glN\>,\glM)$ duality.

\section{The $(\glN\>,\glM)$ duality for KZ and dynamical Hamiltonians}
\label{sec duality}

\subsection{The $(\glN\>,\glM)$ duality}\label{subsec duality}
The Lie algebra $\glN$ acts on $\C[x_1,\dots, x_N]$ by differential operators
$\displaystyle E_{ij}\mapsto x_i\frac \partial{\partial x_j}$. Denote this
$\glN$-module by $\Vb$.
Then
$$
\Vb\,=\,\bigoplus_{m=0}^\infty\,L_{m}^{\langle N\rangle}\ ,
$$
the submodule $L_{m}^{\langle N\rangle}$
being spanned by homogeneous polynomials of degree $m$.
The
$\glN$-module $L_{m}^{\langle N\rangle}$
is irreducible with highest weight $(m,0,\dots,0)$ and
highest weight vector $x_1^{m}$.

We consider $\glN$ and $\glM$
simultaneously. To distinguish generators, modules, etc., we
indicated the dependence on $N$ and $M$ explicitly, for example, $E_{ij}\NNN\>$,
$L_m\MMM $.

Let
$P_{MN} = \C [x_{11},\dots,x_{M1},\dots,x_{1N},\dots,x_{MN}]$
be the space of polynomials of $MN$ variables.
We define the $\glM$-action on $P_{MN}$ by
$\displaystyle E_{ab}\MMM\,
\mapsto \sum_{i=1}^N\,x_{ai} \frac \partial{\partial x_{bi}}$
and the $\glN$-action by
$\displaystyle E_{ij}\NNN\,
\mapsto \sum_{a=1}^M\,x_{ai} \frac \partial{\partial x_{aj}}.
$
There are two isomorphisms of vector spaces,
\bean\label{Miso}
{}&&
\\
\bigl(\C[x_1, \dots, x_M]\big)^{\ox N}\!\!\!\to P_{MN},
&&
(p_1\ox\dots\ox p_N) (x_{11}, \dots,x_{MN})\mapsto
\prod_{i=1}^N\,p_i(x_{1i},\dots, x_{Mi}),
\notag
\\
\bigl(\C[x_1,\dots, x_N]\bigr)^{\ox M}\!\!\! \to P_{MN} ,
&&
(p_1\ox\dots\ox p_M)(x_{11},\dots,x_{MN}) \mapsto
\prod_{a=1}^M\,p_a(x_{a1},\dots, x_{aN}) .
\notag
\eean
Under these isomorphisms, $P_{MN}$ is isomorphic to $(\Vb\MMM) ^{\ox N}$ as a $\glM$-module
and to $({\Vb}^{\langle N \rangle})^{\ox M}$ as a $\glN$-module.

Fix $\bs n=(n_1,\dots,n_N)\in\Zp^{N}$ and
$\bs m=(m_1,\dots, m_M)\in\Zp^{M}$ with $\sum_{i=1}^N n_i = \sum_{a=1}^M m_a$.
Isomorphisms \Ref{Miso} induce an isomorphism of the weight subspaces,
\bean\label{duality isom}
\LNM\ \simeq \LMN\ .
\eean
Under this isomorphism the KZ and dynamical
Hamiltonians interchange,
\bea
H^{{\frak G}, \langle N\rangle}_{a}(\bs\la,\bs z)\,
= \,G^{{\frak X}, \langle M\rangle}_{a}(\bs z,\bs \la)\,,
\qquad
G^{{\rm G}, \langle N\rangle}_{i}(\bs\la,\bs z)\,=\,
H^{{\rm G}, \langle N\rangle}_{i}(\bs z,\bs\la)\,,
\eea
for $a=1,\dots, M$,\ $i=1,\dots,N$,\ \cite{TV5}.

\medskip

Recall that $\;\bs1=(1,\dots, 1)\in\C^N$, and $\xi^{\langle \bs m\rangle}$
is defined by \Ref{xi}.

Let $\bs t^{\langle \bs n\rangle}$ be a Gaudin admissible
critical point of the master function\\
$\Phi(\,\cdot\,;\bs\la-\bs n;\bs m;\bs z)$. Let $(V,\bs z)$
be the associated space of the $(\bs n,\bs\la-\bs n,\bs m,\bs z)$-type.
Let $(U,\bs\la-\bs n+\bs 1)$ be its special bispectral dual space of
the $(\bs m,\bs z,\bs n,\bs\la-\bs n+\bs 1)$-type. Assume that
$(U,\bs\la-\bs n+\bs1)$ is XXX admissible.
Let $\bs t^{\langle \bs m\rangle}+\xi^{\langle \bs m\rangle}$ be the associated
XXX admissible solution of the Bethe ansatz equations \Ref{BAE discrete}.

\medskip

\begin{conj}\label{conjecture}
The corresponding Bethe vectors
$\ \omega^{\frak G}(\bs t^{\langle \bs n\rangle}, \bs z)\in\LMN\;$ and\\
$\omega^{\frak X}(\bs t^{\langle \bs m\rangle},\bs\la)\in\LNM\;$ are
proportional under the duality isomorphism
\Ref{duality isom}.
\hfill $\square$
\end{conj}

\section{The case of $N = M = 2$}
\label{sec m=n}

\subsection{ The $(\glt , \glt)$ duality}
\label { 22 duality}

Let $\bs n=(n_1,n_2)$ and $\bs m=(m_1,m_2)$ be two vectors of nonnegative integers such that
$n_1+n_2=m_1+m_2$.

Fix $\bs \la = (\la_1,\la_2)$ and $\bs z=(z_1,z_2)$ each
with distinct coordinates.

For $m\in\Zp$, let $L_m$ be the irreducible $\glt$-module with highest weight
$(m,0)$ and highest weight vector $v_m$. The vectors $E_{21}^iv_m,
\ i=0,\dots, m$, form a basis in $L_m$.

Set
$$
\alpha = \max\,(0, n_2-m_1)\ , \qquad
\beta =\min \,(m_2,n_2)\ .
$$
The vectors
\bea
\frac{E_{21}^{n_2-i}v_{m_1}}{(n_2-i)!}\ox \frac{E_{21}^iv_{m_2}}{i!}\ ,
\qquad
\alpha\leq i\leq \beta\ ,
\eea
form a basis in
$\LMN\ =\ (L_{m_1}\ox L_{m_2})[n_1,n_2]$.
The vectors
$$
\frac{E_{21}^{m_2-i}v_{n_1}}{(m_2-i)!}\ox \frac{E_{21}^iv_{n_2}}{i!}\ ,
\qquad
\alpha \leq i \leq \beta\ ,
$$
form a basis in
$\LNM\ =\ (L_{n_1}\ox L_{n_2})[m_1,m_2]$.
Isomorphism \Ref{duality isom} identifies the vectors with the same index $i$.

We consider the Gaudin Hamiltonians on $\bs L_{\bs m}[\bs n]$ and the
XXX Hamiltonians on $\bs L_{\bs n}[\bs m]$.

Let $B_{i,j}(t)$ be the series defined in Section
\ref{sec Gaudin dynamical Hamiltonians}. Consider the map
\beq
\label{dyn weyl formula}
b(\bs \la)\ :\
(L_{m_1}\otimes L_{m_2})[n_1,n_2]\ \to\ (L_{m_1}\otimes L_{m_2})[n_2,n_1]\
\eeq
defined by the formula
\bea
\frac{E_{21}^{n_2-i}v_{m_1}}{(n_2-i)!}\ox \frac{E_{21}^iv_{m_2}}{i!}
\ \mapsto\ B_{21}(\la_1-\la_2)
\ \frac{E_{21}^{m_1-n_2+i}v_{m_1}}{(m_1-n_2+i)!}\ox
\frac{E_{21}^{m_2-i}v_{m_2}}{(m_2-i)!} \ .
\eea
For a generic
$\bs \la$ this map is an isomorphism.

By \cite{TV4}, the map $b(\bs \la)
$ commutes with the
Gaudin Hamiltonians in the following sense:
\begin{align*}
b(\bs \la)\, H^{\frak G}_a(\la_1,\la_2,\bs z)\ ={} &
\ H^{\frak G}_a(\la_2, \la_1 ,\bs z)\,b(\bs \la) \ ,
\\
b(\bs \la)\, G^{\frak G}_i(\la_1,\la_2,\bs z)\ ={} &
\ G^{\frak G}_i(\la_2, \la_1 ,\bs z)\,b(\bs \la) \
\end{align*}
for $a, i\,=\,1,2$.

The map
\bean
\frac {E_{21}^{m_2-i}v_{n_1}}
{(m_2-i)!}
\otimes
\frac{E_{21}^iv_{n_2}}{i!}
\ \mapsto\ \frac{E_{21}^{n_1-m_2+i}v_{n_1}}{(n_1-m_2+i)!}\otimes
\frac{E_{21}^{n_2-i}v_{n_2}}{(n_2-i)!}
\eean
defines the Weyl isomorphism
\bean
\label{dual weyl formula}
s\ :\ (L_{n_1}\otimes L_{n_2})[m_1,m_2]\ \to\ (L_{n_1}\otimes L_{n_2})[m_2,m_1]\ .
\eean
The Weyl isomorphism commutes with the XXX Hamiltonians,
\begin{align*}
s \, H^{\frak X}_i(\bs z, \bs \la)\ ={} &
\ H^{\frak X}_i(\bs z, \bs \la)\, s\ ,
\\
s \, G^{\frak X}_a(\bs z, \bs \la)\ ={} &
\ G^{\frak X}_a(\bs z, \bs \la)\, s
\end{align*}
for $i, a\,=\,1,2$.

\subsection{The Gaudin universal rational function and the fundamental
differential equation}
\label{gaudin rational function}
Let $\bs n=(n_1,n_2)$ and $\bs m=(m_1,m_2)$ be two vectors of nonnegative
integers such that $n_1+n_2=m_1+m_2$.

Fix $\bs \la = (\la_1,\la_2)$ and $\bs z=(z_1,z_2)$ each
with distinct coordinates.
{\it The Gaudin universal rational
$\LMN$-valued function} is the function
\bea
\omega (t_1, \dots , t_{n_2})\ =\
\sum_i \,{C_i}\
\frac{E_{21}^{n_2-i}v_{m_1}}{(n_2-i)!}\ \ox \frac{E_{21}^iv_{m_2}}{ i!}\ \ ,
\eea
where
\bean\label{univ function}
C_i(t_1,\dots,t_{n_2})\ = \ \Sym_{n_2}\
\prod_{j=1}^{n_2-i}
\frac 1 {t_j - z_1}\
\prod_{j=1}^{i}
\frac 1 {t_{n_2+j-i} - z_2}
\eean
and $\Sym_n \, f(t_1,\dots , t_{n})\ =\ \sum_{\sigma\in\Sigma_{n}}
f(t_{\sigma(1)},\dots , t_{\sigma(n)})$, see \cite{MaV}.

The Gaudin universal rational function is
symmetric with respect to the group $\Sigma_{n_2}$ of permutations of
variables $t_1,\dots,t_{n_2}$.
The $\Sigma_{n_2}$-orbit of a point
$\bs t^{\langle \bs n \rangle}=(t_1,\dots,t_{n_2})$ is represented by
the polynomial $p(x) = (x-t_1)\dots (x-t_{n_2})$.

The polynomials $(x-z_1)^{n_2-i}(x-z_2)^i$, $i=0,\dots, n_2$, form a basis
in the space of polynomials in $x$ of degree not greater than $n_2$.

\begin{lem}\label{lem y and c_i} Let
$p(x) = (x-t_1)\dots (x-t_{n_2})$ be a polynomial. Let the numbers
$C_0,\dots, C_{n_2}$ be given by formula \Ref{univ function}. Then
$$
\frac{(z_1-z_2)^{n_2}\, p(x)}{p(z_1)\,p(z_2)}\ =\ \sum_{i=0}^{n_2}\
(-1)^i\, C_i\ \frac{(x-z_1)^{n_2-i}}{(n_2-i)!}\ \frac{(x-z_2)^i}{i!}\ .
$$
\hfill $\square$
\end{lem}

Let $V\, =\, \langle \,x^{\la_1} p_1(x)\, ,\, x^{\la_2} p_2(x) \,\rangle$.
Let $(V, \bs z)$ be a space of the $(\bs n, \bs \la, \bs m, \bs z)$-type.
The special fundamental differential operator of $(V, \bs z)$ has the form
\begin{align*}
D\ ={}& \ (x-z_1)(x-z_2) (x\p - \la_1-n_1)(x\p - \la_2)
\\
&{}+\;\phi _{11}\,x\,(x-z_1) (x\p - \la_1-n_1)\;+\;
\phi _{12}\,z_2\, (x-z_1) (x\p - \la_2)
\\
&{}+\;\phi _{21}\, x\, (x-z_2) (x\p - \la_1-n_1)\;+\;
\phi _{22}\,z_1\, (x-z_2) (x\p - \la_2)\ ,
\end{align*}
where $\phi_{ij}$ are suitable numbers such that
\begin{alignat}2
\label{nec cond}
& \phi_{21} + \phi_{22} =\,-\,m_1 , &\qquad\quad
& \phi_{11} + \phi_{12} =\,-\,m_2 ,
\\
& \phi_{22} + \phi_{12} =\,-\,n_1 , && \phi_{21} + \phi_{11} =\,-\,n_2 .
\notag
\end{alignat}
The operator $D$ can be written also in the form
\begin{align*}
D\ ={}& \ (x-z_1)(x-z_2) (x\p - \la_1-n_1)(x\p - \la_2-n_2)
\\
&{}+\;\psi _{11}\, x\,(x-z_1) (x\p - \la_1)\;+\;
\psi _{12}\,z_2\, (x-z_1) (x\p - \la_2-n_2)
\\
&{}+\;\psi _{21}\, x\, (x-z_2) (x\p - \la_1)\;+\;
\psi _{22}\,z_1\, (x-z_2) (x\p - \la_2-n_2)\ ,
\end{align*}
where $\psi_{ij}$ are such that
\begin{alignat}2
& \psi_{21} + \psi_{22} =\,-\,m_1 , &\qquad\quad
& \psi_{11} + \psi_{12} =\,-\,m_2 ,
\notag
\\
& \psi_{22} + \psi_{12} =\,-\,n_1 , && \psi_{21} + \psi_{11} =\,-\,n_2 ,
\notag
\end{alignat}
and \;$\psi_{11}(\la_1-\la_2-n_2)=\phi_{11}(\la_1-\la_2+n_1)+m_2n_2$\ .
\medskip

Equation $Df\, =\, 0$ has a solution
$$
f(x)\ =\ x^{\la_2} \, p_2(x)\ =\ x^{\la_2}\, \sum_{i=0}^{n_2}\,c_i\
\frac{(x-z_1)^{n_2-i}}{(n_2-i)!}\ \frac{(x-z_2)^i}{i!}\ .
$$
Due to conditions \Ref{nec cond}, the expression $x^{-\la_2}Df(x)$
is a polynomial of degree at most $n_2+1$ vanishing at $x=0$.
Expanding this polynomial as a linear combination of the polynomials
$x(x-z_1)^{n_2-i}(x-z_2)^i$, $i=0,\dots, n_2$, we obtain that the equation
$Df=0$ is equivalent to the following relations for the coefficients
$c_0,\dots, c_{n_2}$:
\begin{align}
\label{rec eqns}
z_1 i (m_1 &{}-n_2+i)\,c_{i-1}\;+\;z_2 (m_2-i)(n_2-i)\,c_{i+1}
\\
{}+\bigl(&(z_1+z_2)\,i^2 +
(z_2m_1 - z_1m_2 - 2z_2n_2 - (\la_1-\la_2)(z_1-z_2))\,i
\notag
\\
&\hphantom{(z_1+z_2)\,i^2-(z_1m_2}
-(\la_1-\la_2+n_1)(z_1-z_2)\phi_{11}+z_2m_2n_2\bigr)\ c_i\ =\ 0\ ,
\notag
\end{align}
$i=0,\dots, n_2$. Notice that values $i=0, n_2-m_1, m_2, n_2$ are the values of
$i$ for which equation \Ref{rec eqns} does not contain $c_{i-1}$ or $c_{i+1}$.

\begin{lem}
\label{lem Gaudin closed}
Equations
\Ref{rec eqns} with $i$ such that
$\alpha \leq i\leq \beta $ form a closed
system of equations with respect to
$c_j$ such that
$\alpha \leq j \leq \beta$.
\hfill
$\square$
\end{lem}
Equations \Ref{rec eqns} have a symmetry.
Namely, \Ref{rec eqns} does not change if we replace the parameters
$\la_1,\la_2, n_1, n_2,\phi_{11},i$ by $\la_2,\la_1, n_2, n_1,\psi_{12},m_2-i$,
respectively, and after that replace the unknowns $c_j$ by
$(z_1/z_2)^jc_{m_2-j}$.

\subsection{The XXX universal rational function and the fundamental difference
equation}
\label{XXX rational function}

Let $\bs n=(n_1,n_2)$ and $\bs m=(m_1,m_2)$ be two vectors of nonnegative
integers such that $n_1+n_2=m_1+m_2$.

Fix $\bs \la = (\la_1,\la_2)$ and $\bs z=(z_1,z_2)$ each
with distinct coordinates.
{\it The XXX universal rational
$\LNM$-valued function} is the function
$$
\omega (s_1,\dots, s_{m_2},\bs\la)\ =\
\sum_i
\ {C_i}\
\frac{E_{21}^{m_2-i}v_{n_1}}{(m_2-i)!}\ \ox \frac{E_{21}^iv_{n_2}}{ i!}\ \ ,
$$
where
\begin{alignat}2
\label{XXX univ function}
& C_i(s_1,\dots,s_{m_2},\bs\la)\ ={}
\\[10pt]
& {}=\ \Sym_{m_2}\Biggl[
\ \prod_{a=1}^{m_2-i}\frac 1 {s_a - \la_1+n_1}
& \prod_{b=m_2-i+1}^{m_2}\!\frac {s_b-\la_1}{(s_b-\la_1+n_1)\,(s_b-\la_2+n_2)}&
\notag
\\[4pt]
&&{}\times\,\prod_{a=1}^{m_2-i}\prod_{b=m_2-i+1}^{m_2}\!
\frac {s_a - s_b-1}{s_a - s_b}\ & \Biggr] \ ,
\notag
\end{alignat}
see \cite{KBI, TV2}.

The XXX universal rational function is symmetric with respect to the group
$\Sigma_{m_2}$ of permutations of variables $s_1,\dots,s_{m_2}$. The
$\Sigma_{m_2}$-orbit of a point $s^{\langle \bs n \rangle}=(s_1,\dots,s_{m_2})$
is represented by the polynomial $q(u) = (u-s_1)\dots (u-s_{m_2})$.

The polynomials
$\prod_{j=0}^{m_2-i-1}(u-\la_1+n_1-j)\,\prod_{j=0}^{i-1}(u-\la_2+j)\,$,
$i=0,\dots, m_2$, form a basis
in the space of polynomials in $x$ of degree not greater than $m_2$.

\begin{lem}\label{XXX lem y and c_i} Let
$q(u) = (u-s_1)\dots (u-s_{m_2})$ be a polynomial. Let the numbers
$C_0,\dots, C_{m_2}$ be given by formula \Ref{XXX univ function}. Then
\begin{align*}
& \frac{q(u)}{q(\la_1-n_1)\,q(\la_2-n_2)}\;
\prod_{j=0}^{m_2-1}(\la_1-\la_2+n_2-j)\ ={}
\\[6pt]
& {}=\ \sum_{i=0}^{m_2}\ (-1)^i
\ \frac{C_i(s_1,\dots,s_{m_2},\bs\la)}{(m_2-i)!\,i!}\,
\prod_{j=0}^{m_2-i-1}(u-\la_1+j)\;\prod_{j=0}^{i-1}\,(u-\la_2+n_2-j)\ .
\end{align*}
\hfill $\square$
\end{lem}
Let $U\, =\, \langle \,z_1^{u} q_1(u)\, ,\, z_2^{u} q_2(u) \,\rangle$.
Let $(U, \bs \la)$ be a space of the $(\bs m, \bs z, \bs n, \bs \la)$-type.
The special fundamental difference operator of $(U, \bs \la)$ has the form
\begin{align*}
\widehat D\ ={}& \ (u-\la_1-n_1+1)(u-\la_2+1)(\tau-z_1)(\tau-z_2)
\\
&{}+\;\varphi _{11}\,(u-\la_1-n_1+1)\tau(\tau-z_1)\;+\;
\varphi _{12}\,(u-\la_2+1)z_2(\tau-z_1)\;+\;
\\
&{}+\;\varphi _{21}\,(u-\la_1-n_1+1)\tau(\tau-z_2)\;+\;
\varphi _{22}\,(u-\la_2+1)z_1(\tau-z_2)\ ,
\end{align*}
where $\varphi_{ij}$ are suitable numbers such that
\begin{alignat}2
\label{nec cond2}
& \varphi_{21} + \varphi_{22} =\,-\,m_1 , &\qquad\quad
& \varphi_{11} + \varphi_{12} =\,-\,m_2 ,
\\
& \varphi_{22} + \varphi_{12} =\,-\,n_1 , &&
\varphi_{21} + \varphi_{11} =\,-\,n_2 .
\notag
\end{alignat}

Equation $\widehat Df\, =\, 0$ has a solution
$$
f\ =\ z_2^{u} \, q_2(u)\ =\ z_2^{u}\
\sum_{i=1}^{m_2}\ \frac{c_i}{(m_2-i)!\,i!}
\prod_{j=0}^{m_2-i-1}\!(u-\la_1-n_1+j)\;\prod_{j=0}^{i-1}\,(u-\la_2-j)\ .
$$
Due to conditions \Ref{nec cond2}, the expression $z_2^{-u}\widehat Df(u)$
is a polynomial of degree at most $m_2$. Expanding this polynomial as a linear
combination of the polynomials
$$
\prod_{j=0}^{m_2-i-1}\!(u-\la_1-n_1+j+1)\;\prod_{j=0}^{i-1}\,(u-\la_2-j+1)\ ,
\qquad i=0,\dots, m_2\ ,
$$
we obtain that the equation $\widehat Df=0$ is equivalent
to the following relations for the coefficients $c_0,\dots, c_{m_2}$:
\begin{align}
\label{XXX rec eqns}
z_1 i (m_1 &{}-n_2+i)\,c_{i-1}\;+\;z_2 (m_2-i)(n_2-i)\,c_{i+1}
\\
{}+\bigl(&(z_1+z_2)\,i^2 +
(z_2m_1 - z_1m_2 - 2z_2n_2 - (\la_1-\la_2)(z_1-z_2))\,i
\notag
\\
&\hphantom{(z_1+z_2)\,i^2-(z_1m_2}
-(\la_1-\la_2+n_1)(z_1-z_2)\varphi_{11}+z_2m_2n_2\bigr)\ c_i\ =\ 0\ ,
\notag
\end{align}
$i=0,\dots, m_2$. Notice that values $i=0, n_2-m_1, m_2, n_2$ are the values of
$i$ for which equation \Ref{rec eqns} does not contain $c_{i-1}$ or $c_{i+1}$.

\begin{lem}
Equations \Ref{XXX rec eqns} with $i$ such that $\alpha \leq i\leq \beta $
form a closed system of equations with respect to $c_j$ such that
$\alpha \leq j \leq \beta$.
\hfill $\square$
\end{lem}

Equations \Ref{XXX rec eqns} have a symmetry.
Namely, \Ref{XXX rec eqns} does not change if we replace the parameters
$z_1,z_2, m_1, m_2,\varphi_{11},i$ by $z_2,z_1, m_2, m_1,\varphi_{21},n_2-i$,
respectively, and after that replace the unknowns $c_j$ by $c_{n_2-j}$.
Another symmetry of \Ref{XXX rec eqns} is that the equations do not change
if we replace $\la_1, \la_2$ by $\la_1+1, \la_2+1$.

The pair of equations \Ref{rec eqns} and \Ref{XXX rec eqns} have a symmetry.
Namely, equation \Ref{rec eqns} turns into equation \Ref{XXX rec eqns} if
we replace the parameter $\phi_{11}$ there by $\varphi_{11}$.

\subsection{Proof of Conjecture \ref{conjecture} for $N = M = 2$}
\label{proof of conj}

The symmetries of equations \Ref{rec eqns} and \Ref{XXX rec eqns}
imply the following theorem.

\begin{theorem}
\label{bethe thm}
Let $V\, =\, \langle\, x^{\la_1}p_1(x)\, ,\, x^{la_2} p_2(x)\, \rangle$
and
$U\, =\, \langle \, z_1^u q_1(u) \,,\, z_2^u q_2(u)\, \rangle$.
Let $(V, \bs z)$ be a space of the $(\bs n, \bs \la, \bs m, \bs z)$-type
and $(U, \bs\la+\bs 1)$ the special bispectral dual space
of the $(\bs m, \bs z, \bs n, \bs \la+\bs 1)$-type, where
$\;\bs1=(1,1)\;$. Write
\begin{align*}
p_2(x) &{}\ =\ \sum_{i=0}^{n_2}
\ c_i(\bs\la)\ \frac{(x-z_1)^{n_2-i}}{(n_2-i)!}\ \frac{(x-z_2)^i}{i!}\ ,
\\
p_1(x) &{}\ = \ \sum_{i=0}^{n_1}
\ d_i(\bs\la)\ \frac{(x-z_1)^{n_1-i}}{(n_1-i)!}\ \frac{(x-z_2)^i}{i!}\ ,
\\
q_2(u) &{}\ =\ \sum_{i=1}^{m_2}\ \frac{e_i(\bs\la+\bs1)}{(m_2-i)!\,i!}
\prod_{j=0}^{m_2-i-1}\!(u-\la_1-1-n_1+j)\;\prod_{j=0}^{i-1}\,(u-\la_2-1-j)\ ,
\\
q_1(u) &{}\ =\ \sum_{i=1}^{m_1}\ \frac{f_i(\bs\la+\bs1)}{(m_1-i)!\,i!}
\prod_{j=0}^{m_1-i-1}\!(u-\la_1-1-n_1+j)\;\prod_{j=0}^{i-1}\,(u-\la_2-1-j)\ .
\end{align*}
Then there exist nonzero numbers $A(\bs\la), B(\bs\la), C(\bs\la)$ such that
$$
c_i(\bs\la)\ =\ A(\bs\la)\,d_{m_2-i}(\bs\la)\,(z_1/z_2)^i\ =
\ B(\bs\la)\,e_i(\bs\la+\bs1)\ =\ C(\bs\la)\,f_{n_2-i}(\bs\la+\bs1)
$$
for all \;$i$,\ $\alpha\leq i\leq \beta$.
\end{theorem}

Theorems \ref{bethe thm} and \ref{thm RV} yield the next statement.

\begin{cor}
The Bethe vectors
\begin{alignat*}2
& \sum_{i=\alpha}^{\beta}\ {c_i}(\bs\la+\bs n)
\ \frac{E_{21}^{n_2-i}v_{m_1}}{(n_2-i)!}\ \ox \frac{E_{21}^iv_{m_2}}{ i!}\ ,
&& \sum_{i=\alpha}^{\beta}\ {d_i}(\bs\la+\bs n)
\ \frac{E_{21}^{n_1-i}v_{m_1}}{(n_1-i)!}\ \ox \frac{E_{21}^iv_{m_2}}{ i!}\ ,
\\
& \sum_{i=\alpha}^{\beta}\ {e_i}(\bs\la+\bs n)
\ \frac{E_{21}^{m_2-i}v_{n_1}}{(m_2-i)!}\ \ox \frac{E_{21}^iv_{n_2}}{ i!} \ ,
\qquad && \sum_{i=\alpha}^{\beta}\ {f_i}(\bs\la+\bs n)
\ \frac{E_{21}^{m_1-i}v_{n_1}}{(m_1-i)!}\ \ox \frac{E_{21}^iv_{n_2}}{ i!}
\end{alignat*}
are identified up to proportionality
by isomorphisms \Ref{dyn weyl formula} and \Ref{duality isom}.
\hfill $\square$
\end{cor}

The fact that the first and second Bethe vectors
are identified by isomorphism \Ref{dyn weyl formula}
is proved also in \cite{MV3} in a different way.

\subsection{Correspondence of solutions of the Bethe ansatz equations
for $N = M = 2$}
\label{sec crit points}

Consider the following four systems of the Bethe ansatz equations:
\begin{alignat}2
\label{bea n_2}
&\frac{\la_1-\la_2-1}{t_i}\ +
\ \sum_{a=1}^2 \frac{m_a}{t_i-z_a} \ -\sum_{j=1,\ j\neq i}^{n_2}
\frac{2}{t_i-t_{j}}\ =\ 0\ ,
\qquad && i=1,\dots, n_2\ ,
\\
\label{bea n_1}
& \frac{\la_1-\la_2-1}{t_i}\ +
\ \sum_{a=1}^2 \frac{m_a}{t_i-z_a} \ -\sum_{j=1,\ j\neq i}^{n_1}
\frac{2}{t_i-t_{j}}\ =\ 0\ ,
&& i=1,\dots, n_1\ ,
\\
\label{BAE discrete m_2}
& \quad\ \prod_{i=1}^2\,\frac {t_a -\la_i-1} {t_a-\la_i-1-n_i}\;
\prod_{b, \;b \neq a}^{m_2} \frac {t_a - t_{b} - 1} {t_b - t_{b} + 1}\ =
\ \frac{z_2}{z_1}\ , && a=1,\dots,m_2\ ,
\\
\label{BAE discrete m_1}
& \quad\ \prod_{i=1}^2 \frac {t_a -\la_i-1} {t_a-\la_i-1-n_i}\;
\prod_{b, \;b \neq a}^{m_1} \frac {t_a - t_{b} - 1} {t_b - t_{b} + 1}\ =
\ \frac{z_2}{z_1}\ , && a=1,\dots,m_1\ .
\end{alignat}
The first two systems are of the Gaudin type and the last two systems are
of the XXX type, see \Ref{BAE} and \Ref{BAE discrete}.

Set \,\,$d\,=\,\dim\,(L_{m_1}\otimes L_{m_2})[n_1,n_2] =
\dim\,(L_{m_1}\otimes L_{m_2})[n_2,n_1] =
\dim\,(L_{n_1}\otimes L_{n_2})[m_1,m_2] =
\dim\,(L_{n_1}\otimes L_{n_2})[m_2,m_1]$\;.
By \cite{MV4} for generic $\bs \la$ and $\bs z$
each of the first two systems has exactly $d$ orbits
of Gaudin admissible solutions.
By \cite{MV4} for generic $\bs \la$ and $\bs z$
each of the last two systems has exactly $d$ orbits
of XXX admissible solutions.

\begin{theorem}
\label{thm corresp of solutions}
Let $\bs \la$ and $\bs z$ be generic.
Let $(t_1,\dots, t_{n_2})$ be a Gaudin admissible solution of the system
\Ref{bea n_2}. Let
$(V, \bs z) = \langle x^{\la_1} p_1(x) , x^{\la_2} p_2(x) \rangle$
be the space of the $(\bs \la, \bs z, \bs n, \bs m)$-type associated with the
orbit of the critical point.
Let $(U, \bs \la+\bs1)=\langle z_1^u q_1(u) , z_2^u q_2(u) \rangle$,
where $\bs1=(1,1)$\,, be the special bispectral dual space to $(V, \bs z)$.
Then the polynomial
$p_1$ represents the orbit of a Gaudin admissible solution of the system
\Ref{bea n_1}, the polynomial $q_2$
represents the orbit of an XXX admissible solution of the system
\Ref{BAE discrete m_2}, and the polynomial $q_1$
represents the orbit of an XXX admissible solution of the system
\Ref{BAE discrete m_1}.
\hbox{}\hfill$\square$
\end{theorem}

\section{Baker-Akhieser functions and bispectral correspondence}
\label{sec baker}

\subsection{Grassmannian of Gaudin admissible non-degenerate spaces}
\label{gaudin admissible Grassm}
For $\la \in \C$, we call a vector subspace $X_\la \subset \C[u]$
{\it Gaudin admissible at $\la$} if there exists $m\in\Z_{\geq 0}$ such that
$$
(u-\la)(u-\la-1)\dots (u-\la -m+1)\,\C[u] \subset X_\la
$$
and there exists $f\in X_\la$
such that $f(\la)\neq 0$. We call a vector subspace $X \subset
\C[u]$ {\it Gaudin admissible} if
\bea
X\ =\ \bigcap_{i=1}^n\ X_{\la_i}\
\eea
where $n$ is a natural number, $\la_1,\dots,\la_n\in\C$ are such that
$\la_i-\la_j \notin \Z$ if $i\neq j$,
and for any $i$, $X_{\la_i}$ is Gaudin admissible at $\la_i$.

Each $X_{\la_i}$ can be defined by a finite set of linear
equations. Namely, there exist a positive integer $N_i$, a sequence
of integers $0 < n_{i1}<\dots < n_{iN_i}$ and complex numbers
$c_{i,j,a}$, $a=0,\dots, n_{ij}$, such that $c_{i,j,n_{ij}} \neq 0$
for all $i,j$. Then the space $X_{\la_i}$ consists of all polynomials $r \in
\C[u]$ satisfying the equations
\bean\label{adm Gaudin subspace}
\phantom{aaaaa}
\ \sum_{a=0}^{n_{ij}}\ c_{i,j,a}\ r(\la_i+a)
\ = \ 0\
\quad {\rm for}
\quad
j = 1, \dots , N_i\ .
\eean

For a Gaudin admissible subspace $X$, define the complex vector space $V$ as
the space spanned by functions $x^{\la_i} p_{ij}(x)$, \ $i=1,\dots , n$,\
$j=1,\dots , N_i$, where $p_{ij}(x) = \sum_{a=0}^{n_{ij}}\ c_{i,j,a}\, x^a$.
The space $V$ is a space of quasi-polynomials of the type
considered in Section \ref{subs quasi-polynomials}.

On the other hand, having a space $V$ of quasi-polynomials, as in
Section \ref{subs quasi-polynomials},
we can recover a Gaudin admissible subspace
$X\subset \C[u]$ by formula \Ref{adm Gaudin subspace}.

We say that the space $X$ is {\it non-degenerate} if the corresponding space
$V$ is non-degenerate in the sense of Section \ref{subs quasi-polynomials}.

We denote the set of
all Gaudin admissible non-degenerate subspaces by $\Gr^{\frak G}$ and call it
{\it the Grassmannian} of the Gaudin admissible non-degenerate subspaces.

Let $X\subset\C[u]$ be a Gaudin admissible subspace and $V$ the associated
space of quasi-polynomials in $x$. Define the algebra \ $A_X\, = \,\{\, p\in \C[u]\
| \ p(u)\,X \subset X \,\}.$ An equivalent definition is $A_X\, =\, \{\, p\in
\C[u]\ | \ p(x\partial_x)\,V \subset V\, \}.$

\medskip

Let $\bar D_V$ be the monic fundamental differential operator of $V$.
The function
$$
\Psi_X(x,u) \ =
\ \prod_{i=1}^n\,\prod_{j=1}^{N_i}\,\frac{x}{u-\la_i-n_{ij}}
\ \bar D_V\, x^u \,
$$
is called {\it the stationary Baker-Akhieser function
of the Gaudin admissible space $X$}. Introduce the rational
function $\psi_X (x,u)$ by the formula $\Psi_{X} (x,u)\, =\,
\psi_{X} (x,u)\, x^{u}$.
The function $\psi_{X}(x,u)$
expands as a power series in $x^{-1}, u^{-1}$
of the form
$$
\psi_{X}(x,u) = 1 + \sum_{i,j=1}^\infty \,c_{ij}\, x^{-i}\, u^{-j}\ .
$$
It is easy to see that for every $p\in A_X$, there exists a linear
differential operator $L_p(x,\p)$ with rational coefficients
such that $L_p(x,\p) \,\bar D_V\, =\,
\bar D_V \,p(x\p)$.
As a corollary, we conclude that
$$
L_p(x,\p) \,\Psi_{X}(x,u)\ =\ p(u)\, \Psi_{X}(x,u) \ .
$$
For $p_1, p_2 \in A$, the corresponding operators
$L_{p_1}(x,\p)$ and $L_{p_2}(x,\p)$ commute.

\subsection{Grassmannian of XXX admissible non-degenerate spaces}
\label{XXX admissible Grassm}
Let $z \in \C^*$ be a nonzero complex number with fixed argument.
We call a vector subspace $Y_z \subset \C[x]$
{\it XXX admissible at $z$} if there exists $n \in\Z_{\geq 0}$ such that
$(x-z)^n \,\C[x] \subset Y_z$
and there exists $f\in Y_z$
such that $f(z)\neq 0$. We call a vector subspace $Y \subset
\C[x]$ {\it XXX admissible} if
\bea
Y\ =\ \bigcap_{a=1}^m\ Y_{z_a}\
\eea
where $m$ is a natural number, $z_1,\dots,z_m\in\C^*$ are distinct numbers
and for any $a$, the space $Y_{z_a}$ is XXX admissible at $z_a$.

Each $Y_{z_a}$ can be defined by a finite set of linear equations.
Namely, there exist a positive integer
$M_a$, a sequence of integers $0 < m_{a1}<\dots < m_{aM_a}$ and
complex numbers $d_{a,b,i}$, $i=0,\dots, m_{ab}$, such that
$d_{a,b,m_{ab}} \neq 0$ for all $a,b$. Then
the space $Y_{z_a}$ consists of all polynomials $ r \in \C[x]$
satisfying the equations
\bean\label{adm XXX subspace}
\phantom{aaaaa}
\sum_{i=0}^{m_{ab}}\ d_{a,b,i}\ r^{(i)}(z_a)
\ = \ 0\
\quad
{\rm for}
\quad
b = 1, \dots , M_a\ .
\eean

For an XXX admissible subspace $Y$ define the complex vector space $U$ as
the space spanned by functions $z_a^{u}\, q_{ab}(u)$, $a=1,\dots , m$,\
$b=1,\dots , M_a$, where
\bea
z_a^u\,q_{ab}(u)\ = \ \sum_{i=0}^{m_{ab}}\ d_{a,b,i}\,
\frac {d^i}{du^i}\,x^u\,\vert_{x=z_a}\ .
\eea
To determine $x^u$ at $x=z_a$ we use the chosen argument of $z_a$.

The space $U$ is a space of quasi-exponentials in $u$ of the type
considered in Section \ref{subsec quasi-exponentials}.

On the other hand, having a space $U$ of quasi-exponentials, as in
Section \ref{subsec quasi-exponentials},
we can recover an XXX admissible subspace
$Y\subset \C[x]$ by formula \Ref{adm XXX subspace}.

We say that $Y$ is {\it non-degenerate} if the corresponding $U$ is
non-degenerate in the sense of Section \ref{subsection fund difference oper}.

We denote the set of all XXX admissible non-degenerate subspaces by
$\Gr^{\frak X}$ and call it {\it the Grassmannian} of the XXX
admissible non-degenerate subspaces.

\bigskip

Let $Y\subset\C[x]$ be an XXX admissible subspace and $U$ the associated
space of quasi-exponentials. Define the algebra \ $A_Y = \{ p\in \C[x]\
| \ p(x)\,W \subset W \}.$ An equivalent definition is $A_Y = \{ p\in
\C[x]\ | \ p(\tau_u)\,U \subset U \}.$

Let $\bar D_U$ be the monic fundamental difference operator
of $U$.
The function
\bea
\Phi_{Y} (u,x) \ = \ {\prod_{a=1}^m\,(x-z_a)^{-Ma} \,
\bar D_U\, x^u} \,
\eea
is called {\it
the stationary Baker-Akhieser function
of the XXX admissible space $Y$}. Introduce the rational
function $\phi_{Y} (u,x)$ by the formula
$\Phi_{Y} (u,x)\, =\, \phi_{Y} (u,x)\, x^{u}$.
The function $\phi_{Y}(u,x)$ expands as a power series in $u^{-1}, x^{-1}$
of the form
$$
\phi_{Y}(u,x)\ = \ 1 + \sum_{i,j=1}^\infty \,c_{ij}\, u^{-i}\, x^{-j}\ .
$$
It is easy to see that for every $p\in A_Y$, there exists a linear
difference operator $L_p(u,\tau_u)$ with rational coefficients
such that $L_p(u,\tau_u) \,\bar D_V\, =\,
\bar D_V p(\tau_u)$.
As a corollary, we conclude that
$$
L_p(u,\tau_u) \,\Phi_{Y}(u,x)\ =\ p(x)\, \Phi_{Y}(u,x) \ .
$$
For $p_1, p_2 \in A$, the corresponding operators
$L_{p_1}(u,\tau_u)$ and $L_{p_2}(u,\tau_u)$ commute.

\subsection{Bispectral correspondence}
In terms of the Baker-Akhiezer
functions, we introduce a {\it bispectral}
correspondence between
points of the Grassmanian of
Gaudin admissible non-degenerate subspaces and points
of the Grassmanian of XXX admissible non-degenerate subspaces.
We say that a Gaudin admissible non-degenerate space $X$ corresponds to an
XXX admissible non-degenerate space $Y$ if
$\Psi_{X}(x,u)=\Phi_{Y}(u,x)$. This
correspondence is an analog of Wilson's
bispectral correspondence in \cite{W}.

\begin{theorem}
Let $V$ be a non-degenerate space of quasi-polynomials and let $U$ be
the non-degenerate space of quasi-exponentials which are bispectral dual
with respect to the integral transforms
of Section \ref{Sec Int transforms}.
Then the corresponding admissible spaces
$X \in \Gr^{\frak G}$ and $Y\in \Gr^{\frak X}$
are bispectral correspondent.
\hfill
$\square$
\end{theorem}

\bigskip

\end{document}